\documentclass[12pt,final]{amsart}

\usepackage[margin=1in]{geometry}

\usepackage[dvipsnames]{xcolor}
\usepackage[all]{xy}

\usepackage{amssymb, bm}
\usepackage{amsmath}
\newtheorem{lemma}{Lemma}[section]

\newtheorem{theorem}[lemma]{Theorem}

\newtheorem{corollary}[lemma]{Corollary}

\theoremstyle{definition}
\newtheorem{defn}[lemma]{Definition}

\theoremstyle{remark}
\newtheorem{remark}[lemma]{Remark}
\newtheorem{example}[lemma]{Example}

\newcommand{\Z}{\ensuremath{{\mathbb Z}}}

\newcommand{\N}{\ensuremath{{\mathbb N}}}

\newcommand{\R}{\ensuremath{\mathbb R}}
\newcommand{\C}{\ensuremath{\mathbb C}}

\renewcommand{\Pr}{\ensuremath{\mathbb P}}

\newcommand{\A}{\ensuremath{\mathcal{A}}}
\newcommand{\K}{\ensuremath{\mathcal{K}}}

\DeclareMathOperator{\Tw}{Tw}

\author{Tatyana Barron}
\address{T. Barron, Department of Mathematics, 
University of Western Ontario, 
London, Ontario N6A 5B7, Canada }
\email{tatyana.barron@uwo.ca}
\author{Artour Tomberg}
\address{A. Tomberg, Department of Mathematics, 
University of Western Ontario, 
London, Ontario N6A 5B7, Canada }
\email{atomberg@uwo.ca}

\thanks{Research is supported in part 
by the Natural Sciences and Engineering Research Council of Canada.
}

\date{\today}

\title[]{The twistor space of $\R^{4n}$ and Berezin-Toeplitz operators}

\begin{document}
\sloppy

\maketitle

\noindent {\bf Abstract.} A hyperk\"ahler manifold $M$ has a family of induced complex structures indexed by a two-dimensional sphere $S^2 \cong \mathbb{CP}^1$. The twistor space of $M$ is a complex manifold $\Tw(M)$ together with a natural holomorphic projection $\Tw(M) \to \mathbb{CP}^1$, whose fiber over each point of $\mathbb{CP}^1$ is a copy of $M$ with the corresponding induced complex structure.
We remove one point from this sphere (corresponding to one fiber in the twistor space),
and for the case of $M = \R^{4n}$, $n\in\N$, 
equipped with the standard hyperk\"ahler structure, we construct one quantization
that replaces the family of Berezin-Toeplitz quantizations parametrized by $S^2-\{ pt\}$. We provide semiclassical asymptotics
for this quantization.

\


\


\

\section{Introduction}

Berezin-Toeplitz quantization receives a lot of attention in mathematical literature. Foundational ideas were introduced in papers by Berezin and in works of Boutet de Monvel and Guillemin on Toeplitz
operators. Substantial  recent contributions were made by X. Ma, G. Marinescu, L. Polterovich, 
and many others. Typically, the setting involves Berezin-Toeplitz operators on a symplectic manifold, with a choice of compatible almost complex structure. 
By quantization, we will mean a linear map $f\mapsto T_f^{(k)}$, where $f$ is a function on a symplectic manifold in an appropriate function space, $k\in\N$ is 
a quantum parameter, $T_f^{(k)}$ is a linear operator on a Hilbert space, defined in a such a way that constant functions are mapped to multiples of the identity operator and the Poisson bracket of two functions is mapped to the operator which is asymptotic (as $k\to\infty$) to the commutator of the respective operators. 
  
There is a general question of how quantization depends on the choice of the almost complex structure and there has been a fair amount of effort to investigate this issue. On a hyperk\"ahler manifold, one can ask a related but somewhat different question. Recall that a hyperk\"ahler manifold has a distinguished family of complex structures parametrized by $S^2$, and a corresponding family of K\"ahler forms. At most countably many of these K\"ahler structures are algebraic \cite[Proposition 2.2]{verbitsky:95}, but there exist compact hyperk\"ahler manifolds $(M,g,I_1, I_2,I_3)$ such that the K\"ahler manifold $(M,g,I_j)$ is complex projective for $j=1,2,3$. Examples of such $M$ are tori with linear complex structures and certain $K3$. On such $M$, we have three Berezin-Toeplitz quantizations (each for a different symplectic form), and one can ask how to make one quantization out of these three.  
This was pursued in \cite{barron:17,castejon:16},
see also \cite[Ch.5]{barron:18}. 
For a smooth function $f$ on $M$ there are three Berezin-Toeplitz operators $T_f^{k;j}$, $k\in\N$, 
on each of the K\"ahler manifolds $(M,g,I_j)$. The works  \cite{barron:17, castejon:16} addressed various ways of constructing one
linear operator out of these three and the $k\to\infty$ properties of these operators.

Berezin-Toeplitz quantization on noncompact symplectic manifolds is certainly important, from mathematical and physical
point of view. The basic noncompact case of $\C^n$ is worked out in \cite{coburn:92}. Let us consider $\R^{4n}$ with the standard
hyperk\"ahler structure $(g,I,J,K)$, where $g$ is the standard Euclidean metric, and $I$, $J$, $K$ are three linear complex structures
on  $\R^{4n}$ that satisfy the quaternion relations. Of course (as it is the case in general for hypercomplex manifolds)
there is a whole two-dimensional sphere 
of induced 
complex structures on $\R^{4n}$, 

\begin{equation}
\label{2sphere}
S^2 = \{ aI+bJ+cK \ | \ a,b,c\in\R , \ a^2+b^2+c^2=1\}.
\end{equation}
For each fixed $(a,b,c)$, an  appropriate choice of complex coordinates
identifies $(\R^{4n},aI+bJ+cK)$  with $\C^{2n}$, we get Segal-Bargmann spaces with respect to these complex coordinates,  
and a Berezin-Toeplitz quantization (see the discussion after Theorem \ref{mainth}). 
Instead of following the approach of \cite{barron:17, castejon:16}, we pursue the objective of unifying all these
quantizations (not just three). We can consider the twistor space of $\mathbb{R}^{4n}$, which is a complex manifold $\Tw(\mathbb{R}^{4n})$ parametrizing the different complex structures $aI + bJ + cK$ at points of $\mathbb{R}^{4n}$. There is a hermitian metric on $\Tw(\R^{4n})$ induced by the hyperk\"ahler structure on $\R^{4n}$. This hermitian metric is not K\"ahler (its hermitian form $\omega$ is not closed, see section \ref{subsec:twistorhyper}), although it is balanced. 
The twistor space comes equipped with a natural holomorphic projection $\pi : \Tw(\mathbb{R}^{4n}) \to \mathbb{CP}^1$, where, identifying $\mathbb{CP}^1$ with $S^2$ in the usual way, the fibers of $\pi$ are just the complex manifolds $(\mathbb{R}^{4n}, aI + bJ + cK)$. In this way, the twistor space $\mathrm{Tw}(\mathbb{R}^{4n})$ is, intuitively speaking, the union of these $(\mathbb{R}^{4n}, aI + bJ + cK)$. 
The precise definition is in Section \ref{sec:twistor} below.
Therefore it makes sense to attempt to replace a family of quantizations parametrized by points of $S^2$ by one quantization
on $\Tw
(\R^{4n})$. A general idea to use the twistor space has been in the air (see \cite{barron:17} for discussion in the context of Berezin-Toeplitz quantization).
If $\omega$ were K\"ahler, then the K\"ahler quantization on $\Tw (\R^{4n})$ would have provided a desired answer. But, since $\omega$ is not closed, we can not look for a hermitian holomorphic line bundle with curvature $-i\omega$. This paper is an attempt to put together a "Berezin-Toeplitz-like" quantization on $\Tw (\R^{4n})$ (minus one fiber) by other means.  

Poisson bracket-commutator asymptotics for Berezin-Toeplitz quantization on $\R^{4n}$, with the standard complex structure, appeared in Coburn's work \cite{coburn:92}. We use techniques from \cite{coburn:92}, adjusting for the sphere fiber. However simply "varying complex structure in Theorem 2 \cite{coburn:92}" does not resolve the issue we are addressing, as  it only provides a family of quantizations on $\R^{4n}$ and not one single quantization. This is discussed in some detail after the proof of Theorem \ref{mainth}.

For the proofs, two complex charts on  $\C\Pr^1$ would typically be used, each obtained by deleting a point.
For this reason, we restrict our consideration to the the twistor family which is obtained by  removing a single fiber of the projection $\pi : \Tw(\mathbb{R}^{4n}) \to \mathbb{CP}^1$ from the twistor space $\Tw(\mathbb{R}^{4n})$. The resulting manifold has underlying topological structure $\R^{4n} \times \C$. As an observation, its analytic structure is different from that of a degenerate twistor space \cite{verbitsky:15}, whose underlying topological structure is also the Cartesian product of a hyperk\"ahler manifold with $\C$. This is explained in section \ref{degtwistors}. 
  
The main result is in section \ref{sec:asymp}. From the analysis point view, in its essence, it is a statement about Toeplitz operators
on relevant Bergman spaces. Questions about semi-commutator or commutator or composition of operators, and norm estimates, are standard in theory of Toeplitz and Hankel operators. In complex analysis, these questions are typically posed on domains in $\C^n$. The symbol $f$ is a bounded measurable function on a domain $D$ and the Toeplitz operator $T_f$ with symbol $f$ is an operator on an appropriate function space on $D$. 
Here, philosophically, we are doing the analysis "for many complex structures at the same time". We consider the underlying real domain $D_{\R}$ and a family of complex structures on $D_{\R}$ (that includes the original complex structure), look at the fiber bundle over $D_{\R}$ that parametrizes those complex structures and whose total space has a natural complex structure, define an operator $T_f$ on the total space of this fiber bundle, and work with those operators. 
As it often happens in this subject, there is a strong underlying argument why
this technical statement about operators is useful for quantization. In fact, the motivation comes from differential geometry and mathematical physics, as explained above. 

To discuss the context related to physics more broadly, twistors were introduced in groundbreaking work by Sir Roger Penrose.
Twistor quantization appeared, in particular, in insightful papers \cite{penrose:68, penrose:99}. Let us recall some details.  
Penrose twistors are defined on Minkowski space or, more generally, in a curved spacetime. 
The Minkowski space is a $4$-dimensional real vector space, with a metric of signature $(1,3)$, and its twistor space 
is a $4$-dimensional complex vector space. Discussion of quantization involves
commutation relations for operators that correspond to functions on the twistor space, and the correspondence principle.
Penrose's twistor program ideas have been extended to Riemannian setting and widely applied. Specifically, for $\R^4$, see
the paper on the instanton moduli space by Atiyah, Hitchin and Singer \cite{atiyah:78}.

\section{The twistor space}
\label{sec:twistor}
\subsection{The twistor space of a general hyperk\"ahler manifold $M$
}
\label{subsec:twistorhyper} 
We begin by giving the definition of a hyperk\"ahler manifold and its twistor space, and describing the metric structure of the twistor space.
 
\begin{defn}
A \emph{hyperk\"ahler} manifold is a smooth manifold $M$ together with a triple of integrable almost complex structures $I, J, K : TM \to TM$ satisfying
\[
I^2 = J^2 = K^2 = -\mathrm{Id}, \ IJ = -JI = K,
\]
and a Riemannian metric $g$, simultaneously Hermitian with respect to $I, J, K$, and such that the corresponding Hermitian forms $\omega_I, \omega_J, \omega_K$ are closed.
\end{defn}

It's not hard to verify that with this definition, the form $\omega_J + i\omega_K$ is non-degenerate, has type $(2,0)$ with respect to the structure $I$ and is closed, thus making $(M, I)$ into a holomorphic symplectic manifold.

In addition to $I, J, K$, a hyperk\"ahler manifold has many other complex structures. Any linear combination $A = aI + bJ + cK$ with $a, b, c \in \R$ satisfying $a^2 + b^2 + c^2 = 1$ is also an integrable almost complex structure. The metric $g$ is Hermitian with respect to $A$, and the corresponding Hermitian form $\omega_A$ is closed. In this way, we obtain a family of \emph{induced complex structures} on $M$ parametrized by a two-dimensional sphere:
\[
S^2  = \left\{aI + bJ + cK \ |\ a,b,c\in\R, a^2 + b^2 + c^2 = 1 \right\}.
\]
\begin{defn}
The \emph{twistor space} of a hyperk\"ahler manifold $M$ is the product manifold $\Tw(M) = M \times S^2$.
\end{defn}

Viewing $S^2$ as the set of induced complex structures on $M$ as above, the twistor space $\Tw(M)$ parametrizes these structures at points of $M$. It comes equipped with a natural almost complex sructure $\mathcal{I} : T\Tw(M) \to T\Tw(M)$, defined as follows. Identifying $S^2 \cong \mathbb{CP}^1$ via the stereographic projection (see the next subsection), we let $I_{\mathbb{CP}^1} : T\mathbb{CP}^1 \to T\mathbb{CP}^1$ denote the corresponding almost complex structure. At any point $(x, A) \in M \times \mathbb{CP}^1 \cong \Tw(M)$, the tangent space decomposes as $T_{(x, A)}\Tw(M) = T_x M \oplus T_A \mathbb{CP}^1$; we will call vectors in $T_x M$ vertical and vectors in $T_A \mathbb{CP}^1$ horizontal, and similarly for differential forms. We define
\begin{equation} \label{compstr} 
\begin{array}{ccccc}
\mathcal{I} & : & T_x M \oplus T_A \mathbb{CP}^1 & \longrightarrow & T_x M \oplus T_A \mathbb{CP}^1 \\
&& (X, V) & \longmapsto & \left(AX, I_{\mathbb{CP}^1}V\right).
\end{array}
\end{equation}
It's not hard to verify that $\mathcal{I}^2 = -\mathrm{Id}$, so that $\mathcal{I}$ is an almost complex structure on the twistor space $\Tw(M)$. It actually turns out to be integrable (\cite{salamon:82},
see also \cite{kaledin:98}), making $\Tw(M)$ into a complex manifold. Its complex dimension is clearly one more than the complex dimension of $M$ (in any induced complex structure). For the rest of this subsection, the complex dimension of the twistor space $\Tw(M)$ will be denoted by $m$, so that the complex dimension of $M$ is $m-1$.

The twistor space $\Tw(M) \cong M \times \mathbb{CP}^1$ has two natural projections:
\[
\xymatrix{& \Tw(M) \ar[rd]^\pi \ar[dl]_\sigma & \\ M & & \mathbb{CP}^1,}
\]
the second of which is a holomorphic map. It's not hard to see that for any $A \in \mathbb{CP}^1$, the fiber $\pi^{-1}(A)$ is just the manifold $M$ with the induced complex structure $A$. One can thus think of the twistor space $\Tw(M)$ as the collection of K\"ahler manifolds $(M, A)$ lying above the points $A \in \mathbb{CP}^1$ via the map $\pi$. The sections of $\pi$ are called \emph{twistor lines}.

There is a natural Hermitian metric on the twistor space $\Tw(M)$, namely the product of the hyperk\"ahler metric $g$ from $M$ and the Fubini-Study metric $g_{\mathbb{CP}^1}$ from $\mathbb{CP}^1$:
\[
\sigma^*\left(g\right) + \pi^*\left(g_{\mathbb{CP}^1}\right).
\]
We let $\omega$ denote its Hermitian form and look at its decomposition into its vertical and horizontal parts:
\begin{equation} \label{omega_decomp}
\omega = \omega_M + \omega_{\mathbb{CP}^1}.
\end{equation}
Here both $\omega_M$ and $\omega_{\mathbb{CP}^1}$ are 2-forms on $\Tw(M)$, and at any point of $\Tw(M)$, $\omega_M$ is an element of $\Lambda^2 M$ while $\omega_{\mathbb{CP}^1}$ is an element of $\Lambda^2 \mathbb{CP}^1$. $\omega_{\mathbb{CP}^1}$ is just the pullback of the Fubini-Study form from $\mathbb{CP}^1$ via the projection $\pi$, while $\omega_M$ is not the pullback of any form from $M$, but it has the property that its restriction to any fibre $\pi^{-1}(A) = (M, A)$ is just the K\"ahler form $\omega_A$, defined above. Furthermore, we have the following result.
\begin{lemma} \label{volumeform}
In the above notation, the form $\omega_M^{m-1}$ on the twistor space $\Tw(M)$ is the pullback of a volume form $\Omega$ from $M$ via the map $\sigma : \Tw(M) \to M$:
\[
\omega_M^{m-1} = \sigma^*\left(\Omega\right).
\]
\end{lemma}
\proof
As noted above, the restriction of the form $\omega_M$ on $\Tw(M)$ to the fibre $\pi^{-1}(A)$ is $\omega_A$, in other words,
\[
\left.\omega_M\right|_{\pi^{-1}(A)} = \left.\sigma^*\left(\omega_A\right)\right|_{\pi^{-1}(A)},
\]
and similarly,
\[
\left.\omega_M^{m-1}\right|_{\pi^{-1}(A)} = \left.\sigma^*\left(\omega_A^{m-1}\right)\right|_{\pi^{-1}(A)}.
\]
For any $A \in \mathbb{CP}^1$, $\omega_A^{m-1}$ is a volume form on $M$, since $\omega_A$ is non-degenerate, and the complex dimension of $M$ is $m-1$. Thus, if we show that the forms $\omega_A^{m-1}$ for different $A$ are all equal, the result will follow.

Recall that on any orientable manifold, to each Riemannian metric one can canonically associate two volume forms, and choosing an orientation is equivalent to choosing one of these forms. Indeed, a metric on the tangent bundle induces a metric on the cotangent bundle and all its exterior powers, including the top one, which is a (real) line bundle. A Riemannian metric thus determines two unit vectors in each fiber of the bundle of differential forms of top degree, and in this context an orientation is simply a consistent choice of one of these unit vectors at each point of the manifold.

In our setting of a hyperk\"ahler manifold $M$, the hyperk\"ahler metric $g$ thus gives rise to two volume forms. Each induced complex structure $A$ determines an orientation on $M$, which amounts to choosing one of these two volume forms. We denote this choice by $\mathrm{Vol}_{(M,A)}$. By basic Hermitian geometry,
\[
\mathrm{Vol}_{(M,A)} = \frac{1}{(m-1)!}\ \omega_A^{m-1}.
\]
It only remains to observe that as $A \in \mathbb{CP}^1$ is changing, the right hand side is changing continuously. It follows from this that all forms $\mathrm{Vol}_{(M,A)}$ must be the same, and the same goes for the forms $\omega_A^{m-1}$.
$\Box$

The natural product metric with Hermitian form $\omega$ on the twistor space $\Tw(M)$ described above is never K\"ahler. Indeed, the exterior differential $d$ on $\Tw(M)$ also decomposes into horizontal and vertical parts, $d = d_M + d_{\mathbb{CP}^1}$, and we have, using the decomposition \eqref{omega_decomp},
\[
d\omega = d_M \omega_M + d_{\mathbb{CP}^1} \omega_M + d_M \omega_{\mathbb{CP}^1} + d_{\mathbb{CP}^1} \omega_{\mathbb{CP}^1}.
\]
The first term is zero by the hyperk\"ahler condition on $M$, while the last two terms are zero because $\omega_{\mathbb{CP}^1}$ is a pullback of a closed form from $\mathbb{CP}^1$ to $\Tw(M)$. However, the second term will never be zero (Lemma 4.4 in \cite{kaled-verbit}, see also the proof of Theorem 1 in \cite{tomberg:15}), and thus $d\omega \ne 0$.

On the other hand, the metric on $\Tw(M)$ satisfies the weaker condition of being balanced: $d\left(\omega^{m-1}\right) = 0$. This was first shown by Kaledin and Verbitsky \cite{kaled-verbit}. Indeed,
\[
d\left(\omega^{m-1}\right) = d\left((\omega_M + \omega_{\mathbb{CP}^1})^{m-1}\right) = d\left(\omega_M^{m-1}\right) + (m-1)d\left(\omega_M^{m-2} \wedge \omega_{\mathbb{CP}^1}\right)
\]
Here the first term is zero as a consequence of Lemma \ref{volumeform}. For the second term, we have
\[
d\left(\omega_M^{m-2} \wedge \omega_{\mathbb{CP}^1}\right) = d\left(\omega_M^{m-2}\right) \wedge \omega_{\mathbb{CP}^1} = (m-2)\,d\omega_M \wedge \omega_M^{m-3} \wedge \omega_{\mathbb{CP}^1} =
\]
\[
= (m-2)\left(d_M\omega_M + d_{\mathbb{CP}^1}\omega_M\right) \wedge \omega_M^{m-3} \wedge \omega_{\mathbb{CP}^1} = (m-2)\,d_{\mathbb{CP}^1}\omega_M \wedge \omega_M^{m-3} \wedge \omega_{\mathbb{CP}^1}.
\]
In this last wedge product, $d_{\mathbb{CP}^1}\omega_M \in \Lambda^2 M \otimes \Lambda^1 \mathbb{CP}^1$ at any point of $\Tw(M)$, and since it's being wedged with $\omega_{\mathbb{CP}^1} \in \Lambda^2 \mathbb{CP}^1$, the product will be zero by the dimension of $\mathbb{CP}^1$, which shows that $\Tw(M)$ is balanced.

\subsection{Complex structure and complex coordinates on the twistor space of $\R^{4n}$
}
\label{subsec:twistscs}
Let ${\mathbb{H}}$ denote the algebra of quaternions, with the standard basis
$\{ 1,{\mathrm{i}}, {\mathrm j}, {\mathrm k}\}$
and the relations 
${\mathrm i}^2={\mathrm j}^2={\mathrm k}^2={\mathrm i}{\mathrm j}{\mathrm k}=-1$.
Let $n$ be a positive integer. 
We will use the isomorphisms
\begin{equation}
  \label{isoms}
\R^{4n}\cong \R^4\otimes \R^n\cong {\mathbb{H}}\otimes \R^n\cong \C^{2n}\oplus \C^{2n}{\mathrm {j}}.
\end{equation}
An element of $\R^{4n}$ will be represented by
$$
{\mathbf {x_1}}+{\mathbf {x_2}}{\mathrm{i}}+{\mathbf {x_3}}{\mathrm{j}}+{\mathbf {x_4}}{\mathrm{k}}=
{\mathbf {z}}+{\mathbf {w}}{\mathrm{j}}
,
$$
where
${\mathbf {x_m}}=\begin{pmatrix}x_m^{(1)}\\ ...\\ x_m^{(n)}\end{pmatrix}$ for $m\in\{ 1,2,3,4\}$, and 
$$
{\mathbf {z}}={\mathbf {x_1}}+{\mathbf {x_2}}{\mathrm{i}}, \ {\mathbf {w}}=
{\mathbf {x_3}}+{\mathbf {x_4}}{\mathrm{i}}.
$$
Let $I$, $J$, $K$ be  the three standard linear complex structures on $\R^{4n}$ 
obtained by ${\mathrm i}$, ${\mathrm j}$, ${\mathrm k}$ acting by left
multiplication on ${\mathbb{H}}$. Let $1_n$ denote the $n\times n$ identity matrix, 
and let $0_n$ denote the $n\times n$ zero matrix. 
Using  the standard basis of $\R^4\otimes \R^n$ and the first of the isomorphisms
(\ref{isoms}), we write  
a vector of $\R^{4n}$ as a column vector 
$$
{\mathbf{x}}=\begin{pmatrix}
  {\mathbf {x_1}} \\
  {\mathbf {x_2}} \\
  {\mathbf {x_3}} \\
  {\mathbf {x_4}}
\end{pmatrix}
$$
and observe that 
the matrices
of $I$, $J$, $K$ are, respectively,
\begin{equation}
\label{matIJK}
\begin{pmatrix}0_n & -1_n & & \\
  1_n & 0_n & & \\
   & & 0_n & -1_n  \\
  & & 1_n & 0_n  
\end{pmatrix}, \
\begin{pmatrix} & & -1_n & 0_n  \\
  & & 0_n & 1_n  \\
   1_n & 0_n & &   \\
  0_n & -1_n & &  
\end{pmatrix}, \ 
\begin{pmatrix}&  & & -1_n \\
   &  & -1_n & \\
   & 1_n & &   \\
1_n  & & &   
\end{pmatrix} 
\end{equation}
(an empty spot indicates that the corresponding matrix entries are zero). Or, we can characterize  $I$ as the endomorphism of the tangent bundle
on $\R^{4n}$ such that
\begin{equation}
  \label{Iaction}
I: \ \frac{\partial}{\partial x_1^{(l)} }\mapsto  \frac{\partial}{\partial x_2^{(l)} }, \
\frac{\partial}{\partial x_2^{(l)} }\mapsto  -\frac{\partial}{\partial x_1^{(l)} }, \ \frac{\partial}{\partial x_3^{(l)} }\mapsto  \frac{\partial}{\partial x_4^{(l)} },
 \ \frac{\partial}{\partial x_4^{(l)} }\mapsto  -\frac{\partial}{\partial x_3^{(l)} }, \ l\in \{1,...,n\}
 \end{equation}
 and
 \begin{equation}
   \label{Jaction}
J: \ \frac{\partial}{\partial x_1^{(l)} }\mapsto  \frac{\partial}{\partial x_3^{(l)} }, \
\frac{\partial}{\partial x_3^{(l)} }\mapsto  -\frac{\partial}{\partial x_1^{(l)} }, \ \frac{\partial}{\partial x_4^{(l)} }\mapsto  \frac{\partial}{\partial x_2^{(l)} },
 \ \frac{\partial}{\partial x_2^{(l)} }\mapsto  -\frac{\partial}{\partial x_4^{(l)} }, \ l\in \{1,...,n\}
 \end{equation}
 \begin{equation}
   \label{Kaction}
K: \ \frac{\partial}{\partial x_1^{(l)} }\mapsto  \frac{\partial}{\partial x_4^{(l)} }, \
\frac{\partial}{\partial x_4^{(l)} }\mapsto  -\frac{\partial}{\partial x_1^{(l)} }, \ \frac{\partial}{\partial x_2^{(l)} }\mapsto  \frac{\partial}{\partial x_3^{(l)} },
 \ \frac{\partial}{\partial x_3^{(l)} }\mapsto  -\frac{\partial}{\partial x_2^{(l)} }, \ l\in \{1,...,n\}.
 \end{equation}

The vector space $\R^{4n}$, equipped with the three linear complex structures $I$, $J$, $K$ and the usual flat Euclidean metric, is a hyperk\"ahler manifold. The induced complex structure $aI + bJ + cK$, viewed as an 
endomorphism of $T(\R^{4n})$, takes
$\frac{\partial}{\partial x_1^{(l)} }$ to $a\frac{\partial}{\partial x_2^{(l)} }+b\frac{\partial}{\partial x_3^{(l)} }+c\frac{\partial}{\partial x_4^{(l)} }$,
for $1\le l\le n$, (see (\ref{Iaction}), (\ref{Jaction}), (\ref{Kaction})), 
and so on.
\begin{remark} The space
  of linear complex structures on $\R^{4n}$ is isomorphic to $GL(4n,\R)/GL(2n,\C)$ \cite[I.2.2.5]{mcduff:98}. Here we denote by 
$GL(2n,\C)$ the subgroup of $GL(4n,\R)$ that consists of the matrices of $I$-linear transformations of $\R^{4n}$ or, equivalently, of the matrices in $GL(4n,\R)$ that commute with $I$. Suppose     
  $a_1I+b_1J+c_1K$ and $a_2I+b_2J+c_2K$ are two linear complex structures such that $(a_1,b_1,c_1)\ne (a_2,b_2,c_2)$.
  There is $A\in GL(4n,\R)$ such that $a_1I+b_1J+c_1K=A(a_2I+b_2J+c_2K)A^{-1}$ \cite[Prop. 2.47]{mcduff:98}.  
  
The complex structures $a_1I+b_1J+c_1K$ and $a_2I+b_2J+c_2K$ represent the same point in $GL(4n,\R)/GL(2n,\C)$ (i.e. there is $A\in GL(2n,\C)$ such that $a_1I+b_1J+c_1K=A(a_2I+b_2J+c_2K)A^{-1}$) if and only if $a_1=a_2$. It is straightforward to verify this using the matrix representations (\ref{matIJK}). 
For example, there is no matrix $A$ in $GL(2n,\C)$ such that $AI=IA$ and $AI=JA$, 
  and it is possible to explicitly find $A\in GL(2n,\C)$ such that $AI=IA$ and $AJ=KA$.  
  \end{remark}
For the case of the hyperk\"ahler manifold $\R^{4n}$, its twistor space $\Tw(\R^{4n})$ has a coordinate description, which we now give. 

The sphere (\ref{2sphere}) is a $2$-dimensional sphere in ${\mathrm{End}} (T\R^{4n})$. It is diffeomorphic to the unit sphere $\mathbb{S}^2=\{ (a,b,c)| \ a,b,c\in\R, a^2+b^2+c^2=1\}$ in $\R^3$ via 
$$
s:\mathbb{S}^2\to S^2
$$
$$
(a,b,c)\mapsto aI+bJ+cK .
$$ 
A standard way to 
cover $\mathbb{S}^2\cong \C\Pr^1$ by two complex charts is to introduce
a complex coordinate $\zeta_1$ on $\mathbb{S}^2-\{ pt\}$, taking the point to be, say, $(0,0,1)$,
using the stereographic projection from this point, and 
the complex coordinate $\zeta_2$  on  $\mathbb{S}^2-\{ (0,0,-1)\}$, which is
$\dfrac{1}{\zeta_1}$ on the intersection of the charts.
The formulas for the  
stereographic projection from $(0,0,1)$ (see e.g. \cite[I.\S 6]{conway:78}) give 
$$
\zeta_1=\dfrac{a+ib}{1-c}, \ 
a=\frac{\zeta_1+\bar{\zeta_1}}{|\zeta_1|^2+1}, \ b=\frac{-i(\zeta_1-\bar{\zeta_1})}{|\zeta_1|^2+1}, \
c=\frac{|\zeta_1|^2-1}{|\zeta_1|^2+1}.
$$
We will use a slightly different convention. Let
\begin{equation}
  \label{zetaformula}
\zeta=\dfrac{-c+ib}{a+1}
\end{equation}
be the complex coordinate
on $\mathbb{S}^2-\{ (-1,0,0)\} \cong\C$. It is the complex coordinate $\dfrac{b+ic}{1+a}$, obtained from the stereographic
projection from $(-1,0,0)$, times $i$.  
We note that
$$
a=\frac{1-|\zeta|^2}{1+|\zeta|^2}, \ 
 b=\frac{-i(\zeta-\bar{\zeta})}{|\zeta|^2+1}, \
c=-\frac{\zeta+\bar{\zeta}}{|\zeta|^2+1}. 
$$
On $\mathbb{S}^2-\{ (1,0,0)\}$ we will use the complex coordinate 
$$
\tilde{\zeta}=-\frac{c+ib}{1-a}.
$$
Over $\{ a\ne \pm 1\}$
$$
{\tilde{\zeta}}=\frac{1}{\zeta}.
$$
Let $z_0$ and $z_1$ be the homogeneous coordinates on $\C\Pr^1$. Then 
the diffeomorphisms
$$
\C\Pr^1-\{ [1:0]\}\to \mathbb{S}^2-\{ (-1,0,0)\}, \ \C\Pr^1-\{ [0:1]\}\to \mathbb{S}^2-\{ (1,0,0)\}
$$
are, respectively, obtained from the maps 
$$
\C\Pr^1-\{ [1:0]\}\to \C , \ [z_0:z_1]\mapsto  \zeta=\frac {z_0}{z_1}
$$
and
$$
\C\Pr^1-\{ [0:1]\}\to \C , \ [z_0:z_1]\mapsto  \tilde{\zeta}=\frac {z_1}{z_0}.
$$
Recall that the matrix group $SU(2)$ consists of the matrices  $\begin{pmatrix}
  \alpha & \beta \\
  -\bar\beta & \bar \alpha
\end{pmatrix}$ with $\alpha,\beta\in\C$ such that $|\alpha |^2+|\beta|^2=1$. It acts on $\C\Pr^1$ by
$$
\begin{pmatrix}
  \alpha & \beta \\
  -\bar\beta & \bar \alpha
\end{pmatrix} \ : \ 
              [z_0:z_1]\mapsto [\alpha z_0+\beta z_1:-\bar\beta z_0+\bar\alpha z_1].
              $$
              This action is transitive.
              The corresponding action on $\mathbb{S}^2$ is 
$$
(a,b,c)\mapsto (a',b',c'),
$$
where
$$
a'=(\alpha\bar\alpha-\beta\bar\beta)a+2Im(\alpha\bar\beta)b+2Re(\alpha\bar\beta)c
$$
$$
b'=-i(\alpha\beta-\bar\alpha\bar\beta)a+Re(\alpha^2+\bar\beta^2)b-Im(\alpha^2+\bar\beta^2)c
$$
$$
c'=-(\alpha\beta+\bar\alpha\bar\beta)a+Im(\alpha^2-\bar\beta^2)b+Re(\alpha^2-\bar\beta^2)c.
$$
The twistor space $\Tw(\R^{4n})$ is covered by two charts 
$\R^{4n}\times (S^2-\{ -I\})$ and 
$\R^{4n}\times (S^2-\{ I\})$.
As in \cite{hitchin:13}, we will use complex coordinates $v_1$, ..., $v_n$, $\xi_1$, ...,$\xi_n$,$\zeta$
on $\R^{4n}\times (S^2-\{ -I\})$, where 
\begin{equation}
\label{coordv}
{\mathbf {v}}=\begin{pmatrix}v_1\\ ...\\ v_n\end{pmatrix}={\mathbf {z}}+\zeta\bar{\mathbf {w}}={\mathbf {x_1}}+i{\mathbf {x_2}}+\zeta({\mathbf {x_3}}-i{\mathbf {x_4}})
\end{equation}
and
\begin{equation}
\label{coordksi}
{\bm {\xi}}=\begin{pmatrix}\xi_1\\ ...\\ \xi_n\end{pmatrix}={\mathbf {w}}-\zeta\bar{\mathbf {z}}=
{\mathbf {x_3}}+i{\mathbf {x_4}}-\zeta({\mathbf {x_1}}-i{\mathbf {x_2}}).
\end{equation}
On the chart $\R^{4n}\times (S^2-\{ I\})$, the complex coordinates are
\begin{equation}
\label{tildecoords}
{\tilde{{\mathbf{v}}}}={\tilde{\zeta}}{\mathbf {z}}+\bar{\mathbf {w}}, \ 
{\tilde{{\bm{\xi}}}}={\tilde{\zeta}}{\mathbf {w}}-\bar{\mathbf {z}}, \ {\tilde{\zeta}}. 
\end{equation}
Over the intersection of the charts  
$$
{\mathbf {\tilde{v}}}=\frac{1}{\zeta}{\mathbf {v}}, \ {\tilde{{\bm {\xi}}}}=\frac{1}{\zeta}{\bm \xi}, \
\tilde{\zeta}=\frac{1}{\zeta}. 
$$
\begin{remark}
  The fact that $({\mathbf{v}},{\bm{\xi}},\zeta)$,   $({\tilde{{\mathbf{v}}}}, {\tilde{{\bm{\xi}}}},{\tilde{\zeta}})$
  are complex coordinates on $\Tw
  (\R^{4n})$, for the complex structure (\ref{compstr}), is well known.
  It is used in \cite{hitchin:13, hitchin:14}. If one wishes to verify
  this explicitly, then one can check the equalities
  that involve the almost complex structures and the differentials of the maps.    
\end{remark}

\subsection{Twistor space with one fiber removed}
\label{subsec:twistsfiberrem}
So, we have a diffeomorphism 
$$
Tw(\R^{4n})-(\R^{4n}\times \{ -I\})\to \C^{2n+1}
$$
$$
({\mathbf{x}},aI+bJ+cK)\mapsto
\begin{pmatrix}
  {{\mathbf{v}}}\\
  {{\bm{\xi}}}\\ {\zeta}
  \end{pmatrix}.
  $$
 In  this paper, we will concentrate our attention on
the twistor space of $\R^{4n}$ with one fiber of the projection $\pi : \Tw(\R^{4n}) \to \mathbb{CP}^1$ 
removed, rather than the whole twistor space.
Choosing this fiber, we selected the point to be removed from $S^2$ to be $-I=s((-1,0,0))$
(which is $[1:0]$ in $\C\Pr^1$). We defined complex coordinates, 
$ {{\mathbf{v}}}$,
  ${{\bm{\xi}}}$,  ${\zeta}$ on $Tw(\R^{4n})-(\R^{4n}\times \{ s((-1,0,0))\})$. 

Let $(a_0,b_0,c_0)$ be a point in $\mathbb{S}^2$ such that $(a_0,b_0,c_0)\ne (-1,0,0)$. 
 We will now explain how to define 
complex coordinates, $({\mathbf{v '}},{\bm{\xi '}},\zeta')$, on 
$\Tw(\R^{4n})-(\R^{4n}\times \{ s(a_0,b_0,c_0)\})$.
The point $(a_0,b_0,c_0)$ corresponds to the point $[\zeta_0:1]\in \C\Pr^1$,
where $\zeta_0=\dfrac{-c_0+ib_0}{a_0+1}$. 
The matrix
\begin{equation}
  \label{su2matrix}
\gamma=
\frac{1}{\sqrt{1+|\zeta_0|^2}}
\begin{pmatrix}
  -e^{-i\psi}\zeta_0 &e^{i\psi}\\
  -e^{-i\psi}& -e^{i\psi}\bar\zeta_0
\end{pmatrix}
\end{equation}
where $\psi$ is an arbitrary real number, 
represents an element of $SU(2)$ that takes $[1:0]$ to $[\zeta_0:1]$.
Also 
$$
\gamma^{-1}=
\frac{1}{\sqrt{1+|\zeta_0|^2}}
\begin{pmatrix}
  -e^{i\psi}\bar\zeta_0& -e^{i\psi}\\
  e^{-i\psi}&-e^{-i\psi}\zeta_0   
\end{pmatrix}
$$
takes $[\zeta_0:1]$ to $[1:0]$ and takes $(a_0,b_0,c_0)$ to $(-1,0,0)$.  

Consider a point $P=({\mathbf{x_1}},{\mathbf{x_2}},{\mathbf{x_3}},{\mathbf{x_4}},s(a,b,c))$ in the twistor space, where ${\mathbf{x_1}}$, ${\mathbf{x_2}}$, ${\mathbf{x_3}}$, ${\mathbf{x_4}}$ are arbitrary and $(a,b,c)\ne (a_0,b_0,c_0)$. We set the complex coordinates 
$({\mathbf{v '}},{\bm{\xi '}},\zeta')$, of $P$ to be 
$$
\zeta'=-e^{-2i\psi}\frac{\bar\zeta_0 (-c+ib)+a+1}{-c+ib-\zeta_0(a+1)}
$$
$$
{\mathbf{v '}}=(\zeta' +e^{2i\psi}\bar\zeta_0){\mathbf {z}}+(\zeta_0 \zeta'-e^{2i\psi})\bar{\mathbf {w}}
$$
$$
{\bm{\xi '}}=(\zeta' +e^{2i\psi}\bar\zeta_0){\mathbf {w}}-(\zeta_0 \zeta'-e^{2i\psi})\bar{\mathbf {z}} 
$$
where, as before,
$$
{\mathbf{z}}={\mathbf{x_1}}+{\mathbf{x_2}}i, \ {\mathbf{w}}={\mathbf{x_3}}+{\mathbf{x_4}}i.
$$
To explain, the complex number $\zeta'$ is $\zeta'=\dfrac{-c'+ib'}{a'+1}$, where $(a',b',c')=\gamma^{-1}(a,b,c)$, and  ${\mathbf{v '}}$ 
and ${\bm{\xi '}}$ are (\ref{coordv}) and (\ref{coordksi}) adjusted to the new choice of the complex coordinate on $\mathbb{S}^2-\{ pt\}$ 
(informally speaking, we should express the old coordinate $\zeta$ in terms of the new coordinate $\zeta'$). In particular, if we set $\zeta_0=0$ and $\psi$ 
to be such that $e^{2i\psi}=-1$, then we get ${\mathbf{v '}}=\tilde{{\mathbf{v }}}$,  
${\bm{\xi '}}=\tilde{{\bm{\xi }}}$, and the equalities above become (\ref{tildecoords}).  

Let us denote the map from the twistor space with one fiber removed to $\C^{2n+1}$ defined above by $f$ and denote the standard complex structure on $\C^{2n+1}$ by $J_0$. 
To verify that ${\mathbf{v '}}$, ${\bm{\xi '}}$, $\zeta'$ are indeed complex coordinates, one can explicitly check the equality 
$$
J_0\circ df=df\circ \Bigl ( (aI+bJ+cK)\oplus I_{\mathbb{CP}^1}  \Bigr ).
$$ 
\subsection{Further considerations}
\label{degtwistors}
In this subsection, we make some observations. 

 Let $g$ be the standard flat Euclidean metric on  $\R^{4n}$. On the hyperk\"ahler manifold $(\R^{4n}, g, I,J,K)$, 
 the three K\"ahler forms $\omega_I$, $\omega_J$, $\omega_K$, defined by
 $$
 \omega_I(X,X')=g(IX,X'), \  \omega_J(X,X')=g(JX,X'), \  \omega_K(X,X')=g(KX,X')
$$
 for  all $X,X'\in T_x\R^{4n}$, $x\in \R^{4n}$, have coordinate descriptions  
 $$
 \omega_I=\sum_{l=1}^n\Bigl ( dx_1^{(l)}\wedge dx_2^{(l)}+dx_3^{(l)}\wedge dx_4^{(l)}\Bigr )
 $$
  $$
 \omega_J=\sum_{l=1}^n\Bigl ( dx_1^{(l)}\wedge dx_3^{(l)}+dx_4^{(l)}\wedge dx_2^{(l)}\Bigr )
 $$
  $$
 \omega_K=\sum_{l=1}^n\Bigl ( dx_1^{(l)}\wedge dx_4^{(l)}+dx_2^{(l)}\wedge dx_3^{(l)}\Bigr ).
 $$
The holomorphic symplectic form $ \omega_J+i \omega_K$ is equal to $\sum_{l=1}^ndz_l\wedge dw_l$.  
 
The $2$-form
$\omega_M$ (where $M=\R^{4n}$) on $\Tw
(\R^{4n})$ that appears in the decomposition (\ref{omega_decomp}) is 
defined as follows: at a point $(x,A)$, $x\in \R^{4n}$, $A\in \C\Pr^1$,
for $(X,V), (X',V')\in T_x\R^{4n}\oplus T_A\C\Pr^1$
$$
\omega_M ((X,V), (X',V'))=g(AX,X').
$$
For example, for $n=1$, in the coordinate conventions used above, 
$$
\omega_M=a \ dx_1\wedge dx_2+b \ dx_1\wedge dx_3+c \ dx_1\wedge dx_4+c \ dx_2\wedge dx_3-b \ dx_2\wedge dx_4+a \ dx_3\wedge dx_4.
$$
Also, let us discuss how the set up in this paper relates to definitions in  \cite{verbitsky:15}.  
The twistor family $\R^{4n}\times \Bigl (S^2-\{ pt\}\Bigr )$, 
with its complex structure obtained from $\Tw
(\R^{4n})$, at a quick glance, is possibly reminiscent of the
{\it degenerate
  twistor space} 
of \cite{verbitsky:15}. However, simply removing a point from $S^2$ and obtaining $\C$ instead of $\C\Pr^1$
of complex structures certainly does not have to lead to the situation described in  \cite{verbitsky:15}. 
In \cite{verbitsky:15}, for a compact simple $4n$-dimensional hyperk\"ahler manifold
$(M,g,I,J,K)$,
and a choice of a $(1,1)$-form $\eta$ on $(M,I)$, which is  closed semipositive form of rank $2n$, 
the {\it degenerate
  twistor space} of $M$ is defined as $M\times \C$, with the (integrable) almost complex structure
$I_{\eta}\oplus I_{\C}$, where $I_{\C}$ is the standard complex structure on $\C$, and, at $\zeta\in\C$,
$I_{\eta}$ is the complex structure on $M$ for which $T^{0,1}$ is
\begin{equation}
  \label{condt01}
\{ v\in TM\otimes\C \ | v \lrcorner (\omega_J+i\omega_K+\zeta\eta)=0\} .
\end{equation}
Let us try to explain why our setting  is different from \cite{verbitsky:15}.  Compactness
is not an issue, since linear complex structures descend to the torus $\R^{4n}/\Z^{4n}$. 

Recall that if $A$ is an endomorphism of an even-dimensional vector space $V$ such that $A^2=-1$, 
then $V_{\C}=V\oplus iV$ is $V_{\C}=T^{1,0}\oplus T^{0,1}$, where $T^{1,0}$ is the eigenspace for $i$, 
it consists of the vectors $v-iAv$, $v\in V$, and  $T^{0,1}$ is the eigenspace for $-i$, 
and it consists of the vectors $v+iAv$, $v\in V$. With that in mind, and applying  
(\ref{Iaction}), (\ref{Jaction}), (\ref{Kaction}), we find that 
the complex structure on $\R^{4n}\times \C$  is the one
for which, at a fixed $\zeta\in\C$, 
$T^{0,1}(\R^{4n})$ is the span of $\dfrac{\partial}{\partial \bar z_l}+\zeta\dfrac{\partial}{\partial w_l}$ and
$\dfrac{\partial}{\partial \bar w_l}-\zeta\dfrac{\partial}{\partial z_l}$, $l=1,...,n$. Equivalently, at $\zeta\in\C$, 
$T^{0,1}(\R^{4n})$ is
$$
\{ v\in T(\R^{4n})\otimes \C \ | v\lrcorner \beta=0\}
$$
where
$$
\beta=\sum_{l=1}^n (dz_l+\zeta d\bar w_l)\wedge (dw_l-\zeta d\bar z_l).
$$
Let us try a different choice of the form. Let $n=1$ and choose 
$$
\eta=\frac{i}{2}(dz\wedge d\bar z+dw\wedge d\bar w+dz\wedge d\bar w+dw\wedge d\bar z).
$$
This is a $(1,1)$ form on $(\R^{4n},I)$, which is closed, semipositive in the sense of \cite[Def. 3.6]{verbitsky:15}, 
and of rank $2$ - cf \cite[Def. 3.17]{verbitsky:15}). It defines the complex structure on $\R^4\times \C$ for which
$T^{0,1}$, defined by (\ref{condt01}), is the span of  $\dfrac{\partial}{\partial \bar z}+
\zeta(\dfrac{\partial}{\partial z}-\dfrac{\partial}{\partial w})$ and
$\dfrac{\partial}{\partial \bar w}+
\zeta(\dfrac{\partial}{\partial z}-\dfrac{\partial}{\partial w})$.

\section{Analysis on $\C^{2n+1}$}

In this section, we collect various facts that are necessary to proceed with the proof of the main theorem in section \ref{sec:asymp}.  

We will define several function spaces on $\C^{2n+1}$, $n\in\N$.
 The isomorphism
between $Tw(\R^{4n})-(\R^{4n}\times \{ a_0I+b_0J+c_0K\})$ and $\C^{2n+1}$ described above in section \ref{sec:twistor} will not be used in this section.

Denote by $v_1$,...,$v_n$,$\xi_1$,...$\xi_n$,$\zeta$ the complex coordinates on $\C^{2n+1}$. 
Write, as before, 
$$
{\mathbf {v}}=\begin{pmatrix}v_1\\ ...\\ v_n\end{pmatrix}, \ {\bm {\xi}}=\begin{pmatrix}\xi_1\\ ...\\ \xi_n\end{pmatrix}.
$$
Let $d\mu=d\mu({\mathbf{v}},{\bm{\xi}},\zeta) $ be the Lebesgue measure on $\C^{2n+1}$. We will also write
$$
d\mu( {\mathbf{v}},{\bm{\xi}})=dRe(v_1)dIm(v_1)...dRe(v_n)dIm(v_n)dRe(\xi_1)dIm(\xi_1)...dRe(\xi_n)dIm(\xi_n)
$$
and
$$
d\mu( \zeta)=dRe(\zeta)dIm(\zeta).
$$
We will write,  for two vectors ${\mathbf{a}},{\mathbf{b}}\in\C^n$: 
$$
{\mathbf{a}}\cdot\bar{\mathbf{b}}=\sum_{l=1}^na_l\bar b_l.
$$
\begin{remark}As a comment on notation, for $z\in\C$, $f(z)$ will mean the value of the function $f$ at the point $z$.
  Writing $f(z)$ and not $f(z,{\bar{z}} )$, we are not implying that the function $f$ is holomorphic ($\frac{\partial f}{\partial \bar z}=0$).
  Similarly for $\C^{l}$, $l\in\N$. 
  In the presentation it will be stated which functions are assumed to be holomorphic and it will be clear
  which functions are not assumed to be holomorphic.  
  \end{remark}
Consider the space $L^2(\C^{2n+1},d\mu_k)$, 
where  $k\in\N$ and 
\begin{equation}
\label{measuremuk}
d\mu_k({\mathbf{v}},{\bm{\xi}},\zeta)=(\frac{k}{\pi})^{2n}e^{-k({\mathbf{v}}\cdot\bar{\mathbf{v}}+
  {\bm{\xi}}\cdot\bar{\bm{\xi}})}
  \frac{1}{\pi(1+|\zeta|^2)^2} \ 
d\mu({\mathbf{v}},{\bm{\xi}},\zeta).
\end{equation}
The inner product on $L^2(\C^{2n+1},d\mu_k)$ is
\begin{equation}
  \label{inprodu}
\langle f,g\rangle = \int_{\C^{2n+1}} f({\mathbf{v}},{\bm{\xi}},\zeta)\overline{g({\mathbf{v}},{\bm{\xi}},\zeta)}
d\mu_k({\mathbf{v}},{\bm{\xi}},\zeta)
\end{equation}
and we will write $||.||$ for the corresponding norm. 
Denote by $\A^{(k)}(\C^{2n+1})$ the subspace of $L^2(\C^{2n+1},d\mu_k)$ that consists of holomorphic functions.
  As we show below, these are precisely the 
  holomorphic functions $F$ on $\C^{2n+1}$ that satisfy $\frac{\partial F}{\partial \zeta}\equiv 0$, i.e. 
  depend only on the
  $2n$ complex variables $v_1$,...,$v_n$,$\xi_1$,...,$\xi_n$.  
\begin{lemma}
  \label{lemsupest}
  Suppose $f\in\A ^{(k)}(\C^{2n+1})$. 
  Let $K$ be a compact subset of $\C^{2n+1}$. There is a constant $C_K>0$, depending on $K$ and $n$, such that
  $$
  \underset{z\in K}{\sup}|f(z)|\le C_K ||f||.
  $$
\end{lemma}
\proof
The proof is an obvious modification of the proof of the same statement for the
Bergman space of holomorphic functions in $L^2(\C^{2n+1},d\mu)$ (e.g.  
Lemma 1.4.1 \cite{krantz:92}). We use that the weight
$(\frac{k}{\pi})^{2n}e^{-k\sum_{l=1}^n(|v_l|^2+|\xi_l|^2)}\frac{1}{\pi (1+|\zeta|^2)^2}$ is a positive continuous function on $K$,
therefore it attains its minimum value on $K$, and this value is positive. 
$\Box$
\begin{lemma} The space 
$\A ^{(k)}(\C^{2n+1})$ is a closed subspace of $L^2(\C^{2n+1},d\mu_k)$.
\end{lemma}
\proof
It is sufficient to show that the limit of a sequence of holomorphic functions that converges in $L^2(\C^{2n+1},d\mu_k)$
is a holomorphic function. Suppose $\{ f_n\}$ is such a sequence. We note that $\{ f_n\}$ is Cauchy in norm. 
Pick a compact set $K$ in $\C^{2n+1}$. Restricted to $K$, $\{ f_n\}$ is uniformly Cauchy by Lemma \ref{lemsupest}, and therefore it converges uniformly to a (continuous) function $f_K$. Thus $\{ f_n\}$ converges uniformly on compact sets to a function $f:\C^{2n+1}\to\C$.   
By Theorem 1.9\cite{range:86} (compact convergence of a sequence of holomorphic functions)  $f$ is holomorphic.
$\Box$

Thus, the space $\A ^{(k)}(\C^{2n+1})$ is a Hilbert space with the inner product 
(\ref{inprodu}) (since $L^2(\C^{2n+1},d\mu_k)$ is a Hilbert space, and a closed subspace of a complete metric space is complete). 

Let $l_1$,...,$l_n$,$m_1$,...,$m_n$ be nonnegative integers. Write $l=(l_1,...,l_n)$, $m=(m_1,...,m_n)$. Denote
$$
\psi_{l,m}({\mathbf{v}},{\bm{\xi}},\zeta)=\sqrt{\frac{k^{l_1+...+l_n+m_1+...+m_n}}{l_1!...l_n!m_1!...m_n!}}(v_1)^{l_1}...(v_n)^{l_n}
(\xi_1)^{m_1}...(\xi_n)^{m_n}.
$$
\begin{lemma} The functions $\psi_{l,m}$ 
form a Hilbert space basis (orthonormal basis) in $\A^{(k)}(\C^{2n+1})$.
\end{lemma}
\proof
A holomorphic function on $\C^{2n+1}$ has an expansion 
\begin{equation}
\label{taylorexp}
\sum c_{l_1,...,l_n,m_1,...,m_n,p}(v_1)^{l_1}...(v_n)^{l_n}
(\xi_1)^{m_1}...(\xi_n)^{m_n}\zeta^q,
\end{equation}
where $l_j$, $m_j$, $q$ are nonnegative integers, $c_{l_1,...,l_n,m_1,...,m_n,q}$ are complex numbers (Theorem 1.18 \cite{range:86}, Taylor series of a holomorphic function on a polydisc). 
A monomial in ${\mathbf{v}},{\bm{\xi}}$, times $\zeta^q$ with $q>0$, is not a square integrable function on $\C^{2n+1}$.
More generally, a holomorphic function whose expansion (\ref{taylorexp}) involves a term with  $\zeta^q$ with $q>0$
is not square-integrable. 
It is straightforward to check that the functions $\psi_{l,m}$ form an orthonormal set. 
$\Box$
\begin{lemma}
  For each fixed $p=({\mathbf{v}},{\bm{\xi}},\zeta)\in \C^{2n+1}$, the functional
  $$
  \Phi_{p}:f\mapsto f(p), \ f\in \A^{(k)}(\C^{2n+1})
  $$
  is a continuous linear functional on $\A^{(k)}(\C^{2n+1})$.
\end{lemma}
\proof Take $K$ to be the one point set $\{ p\}$ and apply Lemma \ref{lemsupest}. The statement follows.
$\Box$

Then, by the Riesz representation theorem there is an element $\kappa_{p}\in \A ^{(k)}(\C^{2n+1})$ such that
the linear functional $\Phi_{p}^{(k)}$ is given by the inner product with $\kappa_{p}^{(k)}$:
$$
\langle f,\kappa_{p}^{(k)}\rangle=f(p)
$$
for all $f\in \A^{(k)} (\C^{2n+1})$.  The (weighted) Bergman kernel $\K^{(k)}$ is 
$$
\K^{(k)} (p,q)=\overline{\kappa_{p}^{(k)}(q)}.
$$
It has the reproducing property
\begin{equation}
  \label{reprop1}
h(p)=\int_{\C^{2n+1}} K^{(k)}(p,q)h(q)d\mu_k(q)
\end{equation}
for all $h\in \A ^{(k)}(\C^{2n+1})$. The map
$$
P^{(k)}:L^2(\C^{2n+1},d\mu_k)\to \A ^{(k)}(\C^{2n+1})
$$
$$
f\mapsto \int_{\C^{2n+1}} K^{(k)}(p,q)f(q)d\mu_k(q)
$$
is the orthogonal projection. 
The explicit expression for $\K^{(k)}$
is
$$
\K^{(k)}({\mathbf{u}},{\bm{\eta}},\tau;{\mathbf{v}},{\bm{\xi}},\zeta )=
\sum\psi_{l,m}({\mathbf{u}},{\bm{\eta}})\overline{\psi_{l,m}({\mathbf{v}},{\bm{\xi}})}=
e^{k({\mathbf{u}}\cdot\bar{\mathbf{v}}+
  {\bm{\eta}}\cdot\bar{\bm{\xi}}
  )}.
$$
Thus, the reproducing property (\ref{reprop1}) is
$$
h({\mathbf{u}},{\bm{\eta}},\tau)=\int_{\C^{2n+1}} h({\mathbf{v}},{\bm{\xi}},\zeta )e^{k({\mathbf{u}}\cdot\bar{\mathbf{v}}+
  {\bm{\eta}}\cdot\bar{\bm{\xi}})}d\mu_k({\mathbf{v}},{\bm{\xi}},\zeta)
$$
or
\begin{equation}
  \label{reprop2}
h({\mathbf{u}},{\bm{\eta}},\zeta)=(\frac{k}{\pi})^{2n}
\int_{\C^{2n}} h({\mathbf{v}},{\bm{\xi}},\zeta )e^{k({\mathbf{u}}\cdot\bar{\mathbf{v}}+
  {\bm{\eta}}\cdot\bar{\bm{\xi}})}
e^{-k({\mathbf{v}}\cdot\bar{\mathbf{v}}+
  {\bm{\xi}}\cdot\bar{\bm{\xi}})}
d\mu({\mathbf{v}},{\bm{\xi}}).
\end{equation}
Given $f\in L^2(\C^{2n+1},d\mu_k)$, 
the Toeplitz operator $T_f^{(k)}$, $k\in\N$,  is the linear operator
\begin{equation}
  \label{toeplitzop}
T_f^{(k)}:\A ^{(k)}(\C^{2n+1})\to \A ^{(k)}(\C^{2n+1})
\end{equation}
$$
h\mapsto P^{(k)}(fh).
$$
Usually to $f$ one associates the sequence $\{ T_f^{(k)}\ | k=1,2,3,...\}$. 

So,
$$
(T_f^{(k)}h)(p)=\int_{\C^{2n+1}}f(q)K^{(k)}(p,q)h(q)d\mu_k(q).
$$
\begin{remark}
What we have described just now, is very standard  for the complex euclidean space with  the Gaussian measure (see e.g. \cite{coburn:92}, \cite{berger:86}). 
However, here we consider $\C^{2n+1}$ with the measure (\ref{measuremuk}). 
\end{remark}

The following observation will be useful.  For ${\mathbf{a}}, {\mathbf{b}}\in\C^n$,  the operator defined below is a unitary operator on $\A ^{(k)}(\C^{2n+1})$:
$$
U_{({\mathbf{a}},{\mathbf{b}})}^{(k)}:\A ^{(k)}(\C^{2n+1})\to \A ^{(k)}(\C^{2n+1})
$$
$$
(U_{({\mathbf{a}},{\mathbf{b}})}^{(k)}f)({\mathbf{v}},{\bm{\xi}},\zeta)=
e^{k({\mathbf{v}}\cdot \bar{\mathbf{a}}+{\bm{\xi}}\cdot\bar{\mathbf{b}})
  -\frac{k}{2}({\mathbf{a}}\cdot \bar{\mathbf{a}}+{\mathbf{b}}\cdot \bar{\mathbf{b}})}
f({\mathbf{v-a}},{\bm{\xi}}{\mathbf{-b}},\zeta).
$$

\section{Asymptotic estimates}
\label{sec:asymp}

In this section, we state and prove the main theorem. 

We recall that the twistor space $\Tw (\R^{4n})$ is equipped with two projections: 
\[
\xymatrix{& \Tw(\R^{4n}) \ar[rd]^\pi \ar[dl]_\sigma & \\ \R^{4n} & & \mathbb{CP}^1.}
\]
Choose and fix a point $\tau_0\in \mathbb{CP}^1$. We specified a diffeomorphism
$$
\iota_{\tau_0}:\Tw(\R^{4n})-\pi^{-1}(\tau_0)\to \C^{2n+1}.
$$
Explicitly, if $\tau_0$ corresponds to $(a_0,b_0,c_0)\in\mathbb{S}^2\subset\R^3$ (see sections \ref{subsec:twistscs} and \ref{subsec:twistsfiberrem}), 
we have a map 
\begin{equation}
\label{diffeo}
\R^{4n}\times (\mathbb{S}^2-\{ (a_0,b_0,c_0\})\to \C^{2n+1}.
\end{equation}
If $(a_0,b_0,c_0)=(-1,0,0)$, then this map is given by 
$$
( {\mathbf{x}}, (a,b,c))\mapsto ({\mathbf{v}},{\bm{\xi}},\zeta),
$$
where
\begin{equation}
\label{complexcoord1}
\begin{split}
           \zeta=\frac{-c+ib}{a+1}
         \\
           {\mathbf{v}}={\mathbf{z}}+\zeta\bar{\mathbf{w}}, \ 
           {\bm{\xi}}={\mathbf{w}}-\zeta\bar{\mathbf{z}}
           \end{split}
           \end{equation}
           with the notation
           $$
           {\mathbf{z}}={\mathbf{x_1}}+i{\mathbf{x_2}}, \ {\mathbf{w}}={\mathbf{x_3}}+i{\mathbf{x_4}}.
           $$
           If $(a_0,b_0,c_0)\ne (-1,0,0)$, then let us set $\psi$ in formulas of section  \ref{subsec:twistsfiberrem} so
           that $e^{2i\psi}=-1$, and then 
$$
( {\mathbf{x}}, (a,b,c))\mapsto ({\mathbf{v}},{\bm{\xi}},\zeta),
$$
where
\begin{equation}
\label{complexcoord2}
\begin{split}
\zeta=\frac{\bar\zeta_0 (-c+ib)+a+1}{-c+ib-\zeta_0(a+1)}
\\
{\mathbf{v }}=(\zeta -\bar\zeta_0){\mathbf {z}}+(\zeta_0 \zeta+1)\bar{\mathbf {w}}
\\
{\bm{\xi }}=(\zeta -\bar\zeta_0){\mathbf {w}}-(\zeta_0 \zeta+1)\bar{\mathbf {z}} .
\end{split}
\end{equation}
\begin{remark}
We would like to emphasize that the diffeomorphism $\iota_{\tau_0}$ is \underline{not} an $\R$-linear map between the underlying $(4n+2)$-dimensional real vector spaces.  This is clear from  (\ref{complexcoord1}), (\ref{complexcoord2}). 
\end{remark}
There is a map
\begin{equation}
  \label{mapfunext}
C^{\infty}(\R^{4n})\to C^{\infty}(\C^{2n+1})
\end{equation}
$$
h\mapsto h\circ \sigma\circ (\iota_{\tau_0})^{-1}=\tilde{h}.
$$
Compose this map with the map
\begin{equation}
  \label{funmap2}
  C^{\infty}(\C^{2n+1})\to C^{\infty}(\C^{2n+1})
  \end{equation}
$$
\tilde{h}\to \tilde{h}_{red}
$$
where $\tilde{h}_{red}$ is defined by: at $({\mathbf{v}},{\bm{\xi}},\eta)$  
$$
\tilde{h}_{red}({\mathbf{v}},{\bm{\xi}},\eta)=
\int_{\C}\tilde{h}({\mathbf{v}},{\bm{\xi}},\zeta)\frac{1}{\pi(1+|\zeta|^2)^2}d\mu(\zeta).
$$
It is straightforward to verify that the composition of the two maps above is a well defined linear map $C^{\infty}(\R^{4n})\to C^{\infty}(\C^{2n+1})$ 
(in particular, to conclude that the integral is finite for each ${\mathbf{v}},{\bm{\xi}}$, we can apply 
the Taylor's Theorem (\cite[Th. 2.4.15]{abraham:83} to the function $h$ and use the explicit formulas for ${\mathbf{v}},{\bm{\xi}}$ provided above). 

Let us illustrate (\ref{mapfunext}) and (\ref{funmap2}) by an example. 
\begin{example}
Suppose $n=1$ and the complex coordinates on $\Tw (\R^4)$ with one fiber removed, are, as in (\ref{complexcoord1}) above,   
$$
           v=z+\zeta \bar w, \ 
           \xi=w-\zeta\bar z. 
           $$
Then
$$
x_1+ix_2=z=
\frac{1}{1+\zeta\bar\zeta}(v-\zeta\bar \xi)
$$
$$
x_3+ix_4= w=
  \frac{1}{1+\zeta\bar\zeta}(\xi+\zeta\bar v).
  $$
 Let $h:\R^4\to\C$ be defined by $h(x)=2x_1=z+\bar z$. Then
  $$
  \tilde{h}(v,\xi,\zeta)= \frac{1}{1+\zeta\bar\zeta}(v+\bar v-\zeta\bar\xi-\bar\zeta\xi)
  $$
  and
  $$
  \tilde{h}_{red}(v,\xi,\eta)=\int_{\C}\tilde{h}(v,\xi,\zeta)\frac{1}{\pi(1+|\zeta|^2)^2} \ dRe(\zeta)dIm(\zeta)=4 Re(v)\int_0^{\infty}\frac{rdr}{(1+r^2)^3}
  =Re(v).
  $$
  More generally, a similar calculation for $h(x)=(2x_1)^{2p}$, where $p$ is a positive integer, yields $\tilde{h}_{red}(v,\xi,\eta)=\frac{2^p(Re \ v)^p}{p+1}$.  
\end{example}

For two complex-valued $C^2$ functions on $\C^{2n+1}$ define  
$$
\{ f,g\}  =-i\Bigr (\frac{\partial f}{\partial \zeta}\frac{\partial g}{\partial \bar \zeta}-
\frac{\partial g}{\partial \zeta}\frac{\partial f}{\partial \bar \zeta}+ \sum_{l=1}^n\Bigr (
\frac{\partial f}{\partial v_l}\frac{\partial g}{\partial \bar v_l} -
 \frac{\partial g}{\partial v_l}\frac{\partial f}{\partial \bar v_l} +
 \frac{\partial f}{\partial \xi_l}\frac{\partial g}{\partial \bar \xi_l}-
 \frac{\partial g}{\partial \xi_l}\frac{\partial f}{\partial \bar \xi_l}
\Bigr )\Bigr ).
$$
This bracket is the Poisson bracket for the symplectic form on $\C^{2n+1}$ 
$$
\omega_0=i\Bigr ( d\zeta\wedge d\bar\zeta+\sum_{l=1}^n(dv_l\wedge d\bar v_l+d\xi_l\wedge d\bar\xi_l)\Bigr ) .
$$
On the twistor space with a fiber removed, there is another $2$-form of interest, $\omega_M$: see  (\ref{omega_decomp}) 
and further discussion in section \ref{degtwistors}. Over the fibers of $\pi$, the pull-backs of $\omega_M$ and $\omega_0$ to each fiber are constant multiples of each other (the constant factor depends on the point of the $2$-sphere).

Here is the main result. Assume $f_{\R^{4n}}$ and $g_{\R^{4n}}$ are smooth complex valued functions on $\R^{4n}$, compactly supported. Write $f$ and $g$, respectively, for the functions $\tilde{f}_{\R^{4n}}$ and $\tilde{g}_{\R^{4n}}$ on $\C^{2n+1 }$ defined by (\ref{mapfunext}). 
Then, as we discussed, the corresponding functions $f_{red}$ and $g_{red}$ on $\C^{2n+1}$
defined by (\ref{funmap2})
are well defined. Furthermore, 
\begin{equation}
\label{toeq}  
T_{f}^{(k)}=T_{f_{red}}^{(k)}, \ T_{g}^{(k)}=T_{g_{red}}^{(k)}.
\end{equation}
\begin{theorem}
  \label{mainth}
  There is a constant
  $C=C(f,g)$ such that
  $$
  ||T_f^{(k)}T_g^{(k)}-T_{f_{red}g_{red}}^{(k)}+\frac{1}{k}T^{(k)}_
  {\sum_{j=1}^n\Bigl ( \frac{\partial f_{red}}{\partial v_j}\frac{\partial g_{red}}{\partial \overline{ v _j}} +
      \frac{\partial f_{red}}{\partial \xi_j}\frac{\partial g_{red}}{\partial \overline{ \xi_j}}
      \Bigr )}||
  \le \frac{C}{k^2},
  $$
   $$
  ||T_{f_{red}}^{(k)}T_{g_{red}}^{(k)}-T_{f_{red}g_{red}}^{(k)}+\frac{1}{k}T^{(k)}_
  {\sum_{j=1}^n\Bigl ( \frac{\partial f_{red}}{\partial v_j}\frac{\partial g_{red}}{\partial \overline{ v _j}} +
      \frac{\partial f_{red}}{\partial \xi_j}\frac{\partial g_{red}}{\partial \overline{ \xi_j}}
      \Bigr )}||
  \le \frac{C}{k^2}
  $$
  for all $k>0$.
\end{theorem}
\begin{corollary}
There is a  constant
  $C=C(f,g)$ such that
  $$
  ||ik[T_{f}^{(k)},T_{g}^{(k)}]-T_{\{ f_{red},g_{red}\} }^{(k)}||\le \frac{C}{k}
  $$
  and
   $$
  ||ik[T_{{f_{red}}}^{(k)},T_{{g_{red}}}^{(k)}]-T_{\{ {f_{red}},{g_{red}}\} }^{(k)}||\le \frac{C}{k}
  $$
  for all $k>0$.
\end{corollary}
\noindent   {\bf Proof of Theorem \ref{mainth}.}
    For $\varphi$, $\psi$ in $\A ^{(k)}(\C^{2n+1})$, and writing $\int$ for  $\int_{\C^{2n+1}}$ we have: 
    $$
    \langle T_f^{(k)}T_g^{(k)}\varphi,\psi\rangle=
    \int_{\C^{2n+1}} (T_f^{(k)}T_g^{(k)}\varphi ) ({\mathbf{v}},{\bm{\xi}},\zeta)
    \overline{\psi({\mathbf{v}},{\bm{\xi}},\zeta)}d\mu_k({\mathbf{v}},{\bm{\xi}},\zeta)=
    $$
    $$
    \int \int f({\mathbf{u}},{\bm{\eta}},\tau)
    (T_g^{(k)}\varphi ) ({\mathbf{u}},{\bm{\eta}},\tau)e^{k({\mathbf{v}}\cdot\bar {\mathbf{u}}+
      {\bm{\xi}}\cdot\bar{\bm{\eta}})}d\mu_k({\mathbf{u}},{\bm{\eta}},\tau)
    \overline{\psi({\mathbf{v}},{\bm{\xi}},\zeta)}d\mu_k({\mathbf{v}},{\bm{\xi}},\zeta).
    $$
     Since, by the reproducing property,
    $$
 \int   e^{k({\mathbf{u}}\cdot\bar {\mathbf{v}}+{\bm{\eta}}\cdot\bar{\bm{\xi}})}
    \psi({\mathbf{v}},{\bm{\xi}},\zeta)d\mu_k({\mathbf{v}},{\bm{\xi}},\zeta)= \psi({\mathbf{u}},{\bm{\eta}},\tau)   
    $$
    for each $\tau$, 
    we get:
    $$
    \langle T_f^{(k)}T_g^{(k)}\varphi,\psi\rangle=
 \int  f({\mathbf{u}},{\bm{\eta}},\tau)
    (T_g^{(k)}\varphi ) ({\mathbf{u}},{\bm{\eta}},\tau)
    \overline{\psi({\mathbf{u}},{\bm{\eta}},\tau)}d\mu_k({\mathbf{u}},{\bm{\eta}},\tau)=
    $$
    $$
    \int\int  f({\mathbf{u}},{\bm{\eta}},\tau) g({\mathbf{e}},{\bm{\beta}},\chi)
    \varphi({\mathbf{e}},{\bm{\beta}},\chi)
    e^{k({\mathbf{u}}\cdot\bar {\mathbf{e}}+{\bm{\eta}}\cdot\bar{\bm{\beta}})}d\mu_k({\mathbf{e}},{\bm{\beta}},\chi)
    \overline{\psi({\mathbf{u}},{\bm{\eta}},\tau)}d\mu_k({\mathbf{u}},{\bm{\eta}},\tau).
    $$
    Now let ${\mathbf{a}}={\mathbf{e}}-{\mathbf{u}}$ and ${\mathbf{b}}={\bm{\beta}}-{\bm{\eta}}$. 
    The integral above becomes 
$$
    \int\int  f({\mathbf{u}},{\bm{\eta}},\tau) g({\mathbf{u}}+{\mathbf{a}},{\bm{\eta}}+{\mathbf{b}},\chi)
    \varphi({\mathbf{u}}+{\mathbf{a}},{\bm{\eta}}+{\mathbf{b}},\chi)
    e^{k({\mathbf{u}}\cdot\overline{( {\mathbf{u}}+{\mathbf{a}})}+{\bm{\eta}}\cdot\overline{({\bm{\eta}}+{\mathbf{b}})})}
    $$
    $$
    d\mu_k({\mathbf{u}}+{\mathbf{a}},{\bm{\eta}}+{\mathbf{b}},\chi)
    \overline{\psi({\mathbf{u}},{\bm{\eta}},\tau)}d\mu_k({\mathbf{u}},{\bm{\eta}},\tau)=
    $$
    \begin{equation}
      \label{mainintegral}
    \int\int  f({\mathbf{u}},{\bm{\eta}},\tau) g({\mathbf{u}}+{\mathbf{a}},{\bm{\eta}}+{\mathbf{b}},\chi)
    (U_{-({\mathbf{u}},{\bm{\eta}})}^{(k)}    \varphi)({\mathbf{a}},{\mathbf{b}},\chi)
\end{equation}
    $$
    d\mu_k({\mathbf{a}},{\mathbf{b}},\chi)
    e^{\frac{k}{2}({\mathbf{u}}\cdot\bar{\mathbf{u}}+{\bm{\eta}}\cdot\bar{\bm{\eta}})} 
    \overline{\psi({\mathbf{u}},{\bm{\eta}},\tau)}d\mu_k({\mathbf{u}},{\bm{\eta}},\tau).
    $$
    Now, for $m=2n+3$, we write $g({\mathbf{u}}+{\mathbf{a}},{\bm{\eta}}+{\mathbf{b}},\chi)$ as the sum of the $m$-th Taylor polynomial
    at $({\mathbf{u}},{\bm{\eta}},\chi)$, plus remainder, $g_{m+1}$:
    $$
    g=g-g_{m+1}+g_{m+1}. 
    $$
    First, we aim to establish that
    the part of the integral with $g_{m+1}$ is bounded by $\dfrac{1}{k^2}$ times a positive constant.
    Use the Taylor's Theorem (\cite[Th. 2.4.15]{abraham:83}, Lemma 5 \cite{coburn:92}).  We get: 
  $$
   | \int\int  f({\mathbf{u}},{\bm{\eta}},\tau) g_{m+1}
   (U_{-({\mathbf{u}},{\bm{\eta}})}^{(k)}    \varphi)({\mathbf{a}},{\mathbf{b}},\chi)
    d\mu_k({\mathbf{a}},{\mathbf{b}},\chi)
    e^{\frac{k}{2}({\mathbf{u}}\cdot\bar{\mathbf{u}}+{\bm{\eta}}\cdot\bar{\bm{\eta}})} 
    \overline{\psi({\mathbf{u}},{\bm{\eta}},\tau)}d\mu_k({\mathbf{u}},{\bm{\eta}},\tau)|\le
    $$
    $$
    \sum_{p_1+...+p_{2n}+l_1+...l_{2n}+j_1+j_2=m+1}|c(p_1,...,p_{2n},l_1,...,l_{2n},j_1,j_2)|
    |\partial_{v_1}^{p_1}...\partial_{v_n}^{p_n}\partial_{\xi_1}^{p_{n+1}}...\partial_{\xi_n}^{p_{2n}}
    \partial_{\bar v_1}^{l_1}...\partial_{\bar v_n}^{l_n}\partial_{\bar \xi_1}^{l_{n+1}}...\partial_{\bar \xi_n}^{l_{2n}}
    \partial_{\zeta}^{j_1}\partial_{\bar\zeta}^{j_2}g|_{\infty}
    $$
    $$
    \int(|a_1|^2+...+|a_n|^2+|b_1|^2+...+|b_n|^2)^{\frac{m+1}{2}} 
    |(U_{-({\mathbf{u}},{\bm{\eta}})}^{(k)}    \varphi)({\mathbf{a}},{\mathbf{b}},\chi)|
    d\mu_k({\mathbf{a}},{\mathbf{b}},\chi)
    $$
    $$ \int 
    |f({\mathbf{u}},{\bm{\eta}},\tau)|e^{\frac{k}{2}({\mathbf{u}}\cdot\bar{\mathbf{u}}+{\bm{\eta}}\cdot\bar{\bm{\eta}})} 
    |\psi({\mathbf{u}},{\bm{\eta}},\tau)|d\mu_k({\mathbf{u}},{\bm{\eta}},\tau).
    $$
    Using the Cauchy-Schwarz inequality on the first integral and then the fact that $U_{-({\mathbf{u}},{\bm{\eta}})}^{(k)}$
    is a unitary operator, we conclude that the expression above does not exceed
 $$
    \sum_{p_1+...+p_{2n}+l_1+...l_{2n}+j_1+j_2=m+1}|c(p_1,...,p_{2n},l_1,...,l_{2n},j_1,j_2)|
    |\partial_{v_1}^{p_1}...\partial_{v_n}^{p_n}\partial_{\xi_1}^{p_{n+1}}...\partial_{\xi_n}^{p_{2n}}
    \partial_{\bar v_1}^{l_1}...\partial_{\bar v_n}^{l_n}\partial_{\bar \xi_1}^{l_{n+1}}...\partial_{\bar \xi_n}^{l_{2n}}
    \partial_{\zeta}^{j_1}\partial_{\bar\zeta}^{j_2}g|_{\infty}
    $$
    $$
    ||\varphi|| \ 
    \Bigr (\int(|a_1|^2+...+|a_n|^2+|b_1|^2+...+|b_n|^2)^{m+1}
    d\mu_k({\mathbf{a}},{\mathbf{b}},\chi)\Bigr )^{\frac{1}{2}}
     $$
    $$ \int 
    |f({\mathbf{u}},{\bm{\eta}},\tau)|e^{\frac{k}{2}({\mathbf{u}}\cdot\bar{\mathbf{u}}+{\bm{\eta}}\cdot\bar{\bm{\eta}})} 
    |\psi({\mathbf{u}},{\bm{\eta}},\tau)|d\mu_k({\mathbf{u}},{\bm{\eta}},\tau)\le 
    $$
$$
    \sum_{p_1+...+p_{2n}+l_1+...l_{2n}+j_1+j_2=m+1}|c(p_1,...,p_{2n},l_1,...,l_{2n},j_1,j_2)|
    |\partial_{v_1}^{p_1}...\partial_{v_n}^{p_n}\partial_{\xi_1}^{p_{n+1}}...\partial_{\xi_n}^{p_{2n}}
    \partial_{\bar v_1}^{l_1}...\partial_{\bar v_n}^{l_n}\partial_{\bar \xi_1}^{l_{n+1}}...\partial_{\bar \xi_n}^{l_{2n}}
    \partial_{\zeta}^{j_1}\partial_{\bar\zeta}^{j_2}g|_{\infty}
    $$
    $$
    \Bigr (\int(|a_1|^2+...+|a_n|^2+|b_1|^2+...+|b_n|^2)^{m+1}d\mu_k({\mathbf{a}},{\mathbf{b}},\chi)\Bigr )^{\frac{1}{2}}
     \ ||\varphi||  \ ||\psi||
    $$
    $$
    \Bigl ( \int
    |f({\mathbf{u}},{\bm{\eta}},\tau)|^2 \Bigl (\frac{k}{\pi}\Bigr )^{2n}\frac{1}{\pi(1+|\tau|^2)^2}
    d\mu({\mathbf{u}},{\bm{\eta}},\tau)\Bigr )^{\frac{1}{2}}. 
    $$
    By Lemma 4 \cite{coburn:92}
    $$
    \int(|a_1|^2+...+|a_n|^2+|b_1|^2+...+|b_n|^2)^{m+1}d\mu_k({\mathbf{a}},{\mathbf{b}},\chi)
     =b(m+1,n)\frac{1}{k^{m+1}} 
     $$
     where the constant $b(m+1,n)$ does not depend on $k$. Thus, above we get an estimate with $\dfrac{1}{k^{\frac{m+1}{2}-n}}$,
     which is, for $m=2n+3$, is $\dfrac{1}{k^2}$. This concludes the part of the argument with Taylor remainder.

     Now we consider the part of (\ref{mainintegral}) that involves the $m$-th Taylor polynomial of $g$:
     \begin{equation}
       \label{integralwithpoly}
 \int\int  f({\mathbf{u}},{\bm{\eta}},\tau) (g-g_{m+1})
    (U_{-({\mathbf{u}},{\bm{\eta}})}^{(k)}    \varphi)({\mathbf{a}},{\mathbf{b}},\chi)
    d\mu_k({\mathbf{a}},{\mathbf{b}},\chi)
    e^{\frac{k}{2}({\mathbf{u}}\cdot\bar{\mathbf{u}}+{\bm{\eta}}\cdot\bar{\bm{\eta}})} 
    \overline{\psi({\mathbf{u}},{\bm{\eta}},\tau)}d\mu_k({\mathbf{u}},{\bm{\eta}},\tau).     
\end{equation}
     A typical term
     of the polynomial $g-g_{m+1}$ is of the form
$$
     c(p_1,...,p_{2n},l_1,...,l_{2n},0,0)
    \partial_{u_1}^{p_1}...\partial_{u_n}^{p_n}\partial_{\eta_1}^{p_{n+1}}...\partial_{\eta_n}^{p_{2n}}
    \partial_{\bar u_1}^{l_1}...\partial_{\bar u_n}^{l_n}\partial_{\bar \eta_1}^{l_{n+1}}...\partial_{\bar \eta_n}^{l_{2n}}
    g({\mathbf{u}},{\bm{\eta}},\chi)
    $$
    $$
    a_1^{p_1}...a_n^{p_n}b_1^{p_{n+1}}...b_n^{p_{2n}}
     \bar a_1^{l_1}...\bar a_n^{l_n}\bar b_1^{l_{n+1}}...\bar b_n^{l_{2n}},
     $$
     where  $c(p_1,...,p_{2n},l_1,...,l_{2n},0,0)$ are constants. 
     By Lemma 3 \cite{coburn:92}
$$
     \int a_1^{p_1}...a_n^{p_n}b_1^{p_{n+1}}...b_n^{p_{2n}}
     \bar a_1^{l_1}...\bar a_n^{l_n}\bar b_1^{l_{n+1}}...\bar b_n^{l_{2n}}
     (U_{-({\mathbf{u}},{\bm{\eta}})}^{(k)}    \varphi)({\mathbf{a}},{\mathbf{b}},\chi)
    d\mu_k({\mathbf{a}},{\mathbf{b}},\chi)=
    $$
    \begin{equation}
      \label{intermintegral}
    \frac{1}{k^{l_1+...+l_{2n}}}\int_{\C} \partial_{a_1}^{l_1}...\partial_{a_n}^{l_n}\partial_{b_1}^{l_{n+1}}...\partial_{b_n}^{l_{2n}} \Bigl [ 
    a_1^{p_1}...a_n^{p_n}b_1^{p_{n+1}}...b_n^{p_{2n}} (U_{-({\mathbf{u}},{\bm{\eta}})}^{(k)}    \varphi)({\mathbf{a}},{\mathbf{b}},\chi) \Bigr ] 
    \Bigr |_{{\mathbf{a}}={\mathbf{b}}={\mathbf{0}}} \ \frac{d\mu(\chi)}{\pi(1+|\chi|^2)^2}.
    \end{equation}
    If there is a $j$ such that $l_j<p_j$, then this expression is equal to zero. Assume $l_j\ge p_j$ for all
    $j\in \{ 1,...,2n\}$. Then (\ref{intermintegral}) is a sum of the terms, with coefficients independent of $k$,
    that are
    $$
     \frac{1}{k^{l_1+...+l_{2n}}}\int_{\C} \partial_{a_1}^{t_1}...\partial_{a_n}^{t_n}\partial_{b_1}^{t_{n+1}}...\partial_{b_n}^{t_{2n}}
     (U_{-({\mathbf{u}},{\bm{\eta}})}^{(k)}    \varphi)({\mathbf{a}},{\mathbf{b}},\chi)
    \Bigr |_{{\mathbf{a}}={\mathbf{b}}={\mathbf{0}}} \ \frac{d\mu(\chi)}{\pi(1+|\chi|^2)^2}
    $$
    where all $0\le t_j\le l_j$, and by Lemma 2 \cite{coburn:92} this is
 $$
    \frac{1}{k^{l_1+...+l_{2n}}}\int_{\C}
e^{\frac{k}{2}({\mathbf{u}}\cdot\bar {\mathbf{u}}+{\bm{\eta}}\cdot\bar {\bm{\eta}})}
    \partial_{u_1}^{t_1}...\partial_{u_n}^{t_n}\partial_{\eta_1}^{t_{n+1}}...\partial_{\eta_n}^{t_{2n}}
     \{  \varphi ({\mathbf{u}},{\bm{\eta}},\chi)e^{-k({\mathbf{u}}\cdot\bar {\mathbf{u}}+{\bm{\eta}}\cdot\bar {\bm{\eta}})}\}
    \frac{d\mu(\chi)}{\pi(1+|\chi|^2)^2}.
    $$
    Putting this in (\ref{integralwithpoly}), we get terms (with coefficients that do not depend on $k$)
    \begin{equation}
      \label{nextintegral}
     \frac{1}{k^{l_1+...+l_{2n}}}\int_{\C}
 \int  f({\mathbf{u}},{\bm{\eta}},\tau) 
 \overline{\psi({\mathbf{u}},{\bm{\eta}},\tau)}
  \partial_{u_1}^{p_1}...\partial_{u_n}^{p_n}\partial_{\eta_1}^{p_{n+1}}...\partial_{\eta_n}^{p_{2n}}
    \partial_{\bar u_1}^{l_1}...\partial_{\bar u_n}^{l_n}\partial_{\bar \eta_1}^{l_{n+1}}...\partial_{\bar \eta_n}^{l_{2n}}
    g({\mathbf{u}},{\bm{\eta}},\chi)
    \end{equation}
    $$
 \partial_{u_1}^{t_1}...\partial_{u_n}^{t_n}\partial_{\eta_1}^{t_{n+1}}...\partial_{\eta_n}^{t_{2n}}
     \{  \varphi ({\mathbf{u}},{\bm{\eta}},\chi)e^{-k({\mathbf{u}}\cdot\bar {\mathbf{u}}+{\bm{\eta}}\cdot\bar {\bm{\eta}})}\}
\frac{1}{\pi(1+|\tau|^2)^2}
(\frac{k}{\pi})^{2n}
     d\mu({\mathbf{u}},{\bm{\eta}},\tau) \frac{d\mu(\chi)}{\pi(1+|\chi|^2)^2}.
     $$
     Recall that $l_j\ge p_j$ and $l_j\ge t_j$ for all $j\in \{ 1,...,2n\}$ and $p_1+...+p_{2n}+l_1+...+l_{2n}\le m$.
     After complex integration by parts (taking into account that $f$ and $g$ have compact support
     with respect to the variables $u_1$, ... ,$u_n$, $\eta_1$, ..., $\eta_n$), we conclude that 
     (\ref{nextintegral}) becomes a sum of the terms
     (with coefficients that do not depend on $k$) 
     $$
  \frac{1}{k^{l_1+...+l_{2n}}}\int_{\C} \int    
  \varphi ({\mathbf{u}},{\bm{\eta}},\chi)\overline{\psi({\mathbf{u}},{\bm{\eta}},\tau)}
   \partial_{u_1}^{p_1+s_1}...\partial_{u_n}^{p_n+s_n}\partial_{\eta_1}^{p_{n+1}+s_{n+1}}...\partial_{\eta_n}^{p_{2n}+s_{2n}}
    \partial_{\bar u_1}^{l_1}...\partial_{\bar u_n}^{l_n}\partial_{\bar \eta_1}^{l_{n+1}}...\partial_{\bar \eta_n}^{l_{2n}}
    g({\mathbf{u}},{\bm{\eta}},\chi)
    $$
    $$
 \partial_{u_1}^{q_1}...\partial_{u_n}^{q_n}\partial_{\eta_1}^{q_{n+1}}...\partial_{\eta_n}^{q_{2n}} f({\mathbf{u}},{\bm{\eta}},\tau)
     d\mu_k({\mathbf{u}},{\bm{\eta}},\tau) \frac{d\mu(\chi)}{\pi(1+|\chi|^2)^2}.
     $$
     Thus, the terms with $l_1+...+l_{2n}>1$ lead to the inequality $\le \dfrac{const}{k^2}$. It remains to consider
     the cases $l_1+...+l_{2n}=0$ and $l_1+...+l_{2n}=1$.

     Looking at (\ref{intermintegral}), we conclude that the part of (\ref{integralwithpoly})  with $l_1+...+l_{2n}=0$ is
     \begin{equation}
       \label{zeroterm}
  \int\int  f({\mathbf{u}},{\bm{\eta}},\tau) g({\mathbf{u}},{\bm{\eta}},\chi) 
    (U_{-({\mathbf{u}},{\bm{\eta}})}^{(k)}    \varphi)({\mathbf{a}},{\mathbf{b}},\chi)
    d\mu_k({\mathbf{a}},{\mathbf{b}},\chi)
    e^{\frac{k}{2}({\mathbf{u}}\cdot\bar{\mathbf{u}}+{\bm{\eta}}\cdot\bar{\bm{\eta}})} 
    \overline{\psi({\mathbf{u}},{\bm{\eta}},\tau)}d\mu_k({\mathbf{u}},{\bm{\eta}},\tau). 
    \end{equation}
    Setting ${\mathbf{u}}+{\mathbf{a}}={\mathbf{v}}$ and  ${\bm{\eta}}+{\mathbf{b}}={\bm{\xi}}$, and using
    the reproducing property (\ref{reprop2}), we observe that for each ${\mathbf{u}}$, ${\bm{\eta}}$ 
$$
    \int  g({\mathbf{u}},{\bm{\eta}},\chi) 
 (U_{-({\mathbf{u}},{\bm{\eta}})}^{(k)}    \varphi)({\mathbf{a}},{\mathbf{b}},\chi)
    d\mu_k({\mathbf{a}},{\mathbf{b}},\chi) e^{\frac{k}{2}({\mathbf{u}}\cdot\bar{\mathbf{u}}+{\bm{\eta}}\cdot\bar{\bm{\eta}})} =
    \int_{\C}  g({\mathbf{u}},{\bm{\eta}},\chi)  \varphi({\mathbf{u}},{\bm{\eta}},\chi)\frac{d\mu(\chi)}{\pi(1+|\chi|^2)^2}.
    $$
So, we rewrite (\ref{zeroterm}) as 
    $$
    \int\int_{\C}  f({\mathbf{u}},{\bm{\eta}},\tau) g({\mathbf{u}},{\bm{\eta}},\chi) 
       \varphi({\mathbf{u}},{\bm{\eta}},\chi)\frac{d\mu(\chi)}{\pi(1+|\chi|^2)^2}
    \overline{\psi({\mathbf{u}},{\bm{\eta}},\tau)}d\mu_k({\mathbf{u}},{\bm{\eta}},\tau). 
    $$
    Now suppose $l_1+...+l_{2n}=1$. Without loss of generality $l_1=1$ and $l_2=...=l_{2n}=0$. Then $p_1=0,1$ and 
    $p_2=...=p_{2n}=0$. The corresponding integral is
         \begin{equation}
       \label{firstorderterms}
       \int\int  f({\mathbf{u}},{\bm{\eta}},\tau) \Bigl ( \partial _{\bar u_1} g({\mathbf{u}},{\bm{\eta}},\chi)\bar a_1+
     \partial _{u_1}  \partial _{\bar u_1}g({\mathbf{u}},{\bm{\eta}},\chi)a_1\bar a_1\Bigr ) 
    (U_{-({\mathbf{u}},{\bm{\eta}})}^{(k)}    \varphi)({\mathbf{a}},{\mathbf{b}},\chi)
     d\mu_k({\mathbf{a}},{\mathbf{b}},\chi)
 \end{equation}
     $$
    e^{\frac{k}{2}({\mathbf{u}}\cdot\bar{\mathbf{u}}+{\bm{\eta}}\cdot\bar{\bm{\eta}})} 
    \overline{\psi({\mathbf{u}},{\bm{\eta}},\tau)}d\mu_k({\mathbf{u}},{\bm{\eta}},\tau). 
    $$
    By Lemma 3 \cite{coburn:92}
     $$
    \int (U_{-({\mathbf{u}},{\bm{\eta}})}^{(k)}    \varphi)({\mathbf{a}},{\mathbf{b}},\chi)a_1\bar a_1 \  d\mu_k({\mathbf{a}},{\mathbf{b}},\chi)=
    \frac{1}{k}\int_{\C}\partial _{a_1}(a_1U_{-({\mathbf{u}},{\bm{\eta}})}^{(k)}    \varphi)({\mathbf{a}},{\mathbf{b}},\chi)
    \Bigr |_{{\mathbf{a}}={\mathbf{b}}={\mathbf{0}}} \ \frac{d\mu(\chi)}{\pi(1+|\chi|^2)^2}=
    $$
    $$
    \frac{1}{k}\int_{\C}(U_{-({\mathbf{u}},{\bm{\eta}})}^{(k)}    \varphi)({\mathbf{0}},{\mathbf{0}},\chi)
    \ \frac{d\mu(\chi)}{\pi(1+|\chi|^2)^2}=
    \frac{1}{k}\int_{\C}
  e^{-\frac{k}{2}({\mathbf{u}}\cdot\bar{\mathbf{u}}+{\bm{\eta}}\cdot\bar{\bm{\eta}})}
    \varphi ({\mathbf{u}},{\bm{\eta}},\chi)
    \ \frac{d\mu(\chi)}{\pi(1+|\chi|^2)^2}
     $$
     and
         $$
    \int (U_{-({\mathbf{u}},{\bm{\eta}})}^{(k)}    \varphi)({\mathbf{a}},{\mathbf{b}},\chi)\bar a_1  \  d\mu_k({\mathbf{a}},{\mathbf{b}},\chi)=
    \frac{1}{k}\int_{\C}\partial _{a_1}(U_{-({\mathbf{u}},{\bm{\eta}})}^{(k)}    \varphi)({\mathbf{a}},{\mathbf{b}},\chi)
    \Bigr |_{{\mathbf{a}}={\mathbf{b}}={\mathbf{0}}} \ \frac{d\mu(\chi)}{\pi(1+|\chi|^2)^2}
    $$
    which, by Lemma 2 \cite{coburn:92}, becomes
         $$
    \frac{1}{k}\int_{\C}
 e^{\frac{k}{2}({\mathbf{u}}\cdot\bar{\mathbf{u}}+{\bm{\eta}}\cdot\bar{\bm{\eta}})}
    \partial _{u_1}\{ \varphi({\mathbf{u}},{\bm{\eta}},\chi) e^{-k({\mathbf{u}}\cdot\bar{\mathbf{u}}+{\bm{\eta}}\cdot\bar{\bm{\eta}})}\}
     \ \frac{d\mu(\chi)}{\pi(1+|\chi|^2)^2}.
    $$  
     Thus, (\ref{firstorderterms}) is equal to 
     $$
\frac{1}{k}
\int\int_{\C}  f({\mathbf{u}},{\bm{\eta}},\tau) \Bigl ( \partial _{\bar u_1} g({\mathbf{u}},{\bm{\eta}},\chi)
 e^{\frac{k}{2}({\mathbf{u}}\cdot\bar{\mathbf{u}}+{\bm{\eta}}\cdot\bar{\bm{\eta}})}
    \partial _{u_1}\{ \varphi({\mathbf{u}},{\bm{\eta}},\chi) e^{-k({\mathbf{u}}\cdot\bar{\mathbf{u}}+{\bm{\eta}}\cdot\bar{\bm{\eta}})}\} +
$$
$$
    \partial _{u_1}  \partial _{\bar u_1}g({\mathbf{u}},{\bm{\eta}},\chi)
  e^{-\frac{k}{2}({\mathbf{u}}\cdot\bar{\mathbf{u}}+{\bm{\eta}}\cdot\bar{\bm{\eta}})}
    \varphi ({\mathbf{u}},{\bm{\eta}},\chi)
    \Bigr ) 
    e^{\frac{k}{2}({\mathbf{u}}\cdot\bar{\mathbf{u}}+{\bm{\eta}}\cdot\bar{\bm{\eta}})} 
    \overline{\psi({\mathbf{u}},{\bm{\eta}},\tau)}
 \ \frac{d\mu(\chi)}{\pi(1+|\chi|^2)^2}
    d\mu_k({\mathbf{u}},{\bm{\eta}},\tau).  
    $$
    After a complex integration by parts on the first summand (taking into account that $f$ and $g$ are compactly supported
    on $\C^{2n}$), we get
     $$
-\frac{1}{k}
\int\int_{\C}\partial_{u_1}  f({\mathbf{u}},{\bm{\eta}},\tau)  \partial _{\bar u_1} g({\mathbf{u}},{\bm{\eta}},\chi)
 \varphi({\mathbf{u}},{\bm{\eta}},\chi) 
    \overline{\psi({\mathbf{u}},{\bm{\eta}},\tau)}
 \ \frac{d\mu(\chi)}{\pi(1+|\chi|^2)^2}
    d\mu_k({\mathbf{u}},{\bm{\eta}},\tau).  
    $$   
To summarize, we got: for arbitrary $\varphi$, $\psi$ in $\A ^{(k)}(\C^{2n+1})$
$$
|\langle \Bigl (T_f^{(k)}T_g^{(k)} -  T_{f_{red}g_{red}}^{(k)}+\frac{1}{k}T^{(k)}_
  {\sum_{j=1}^n\Bigl ( \frac{\partial f_{red}}{\partial v_j}\frac{\partial g_{red}}{\partial \overline{ v _j}} +
      \frac{\partial f_{red}}{\partial \xi_j}\frac{\partial g_{red}}{\partial \overline{ \xi_j}} \Bigr ) }
       \Bigr )\varphi       ,\psi \rangle  |\le \frac{const}{k^2}.
$$  
We recall that for a bounded selfadjoint operator $A$ on a Hilbert space $H$ we have $||A||=\underset{||x||=1}{\sup}|\langle Ax,x\rangle |$. This equality also holds  if $A$ is bounded and $iA$ is selfadjoint. Then, for an arbitrary bounded operator $B$ on $H$ we write $B=\frac{1}{2}(B+B^*+B-B^*)$ and conclude that  
$||B||\le 2\underset{||x||=1}{\sup}|\langle Bx,x\rangle |$.

This proves the first statement of the theorem. The second statement follows from (\ref{toeq}). 
    $\Box$

\begin{remark}
The proof is a modification of the proof of Theorem
2 in \cite{coburn:92}. This way of obtaining asymptotics for Toeplitz operators is quite common: write an inner product as an integral, break the symbol into 
a Taylor polynomial and remainder, estimate the two parts separately. 
\end{remark}
\begin{remark}
In the proof, it was only used that  ${\mathbf{v}}$, ${\bm{\xi}}$, $\zeta$ are complex coordinates on $\Tw(\R^{4n})$ with one fiber removed. We did not utilize the explicit formulas (\ref{complexcoord1}) or (\ref{complexcoord2}). 
\end{remark}
\begin{remark}
In the assumptions of the theorem, $f$ and $g$ are the functions obtained from $f_{\R^{4n}}$ and $g_{\R^{4n}}$, smooth compactly supported functions on $\R^{4n}$. In fact, for the proof it was only needed that 
$f_{\R^{4n}}\in C_c^{l}(\R^{4n})$ (the set of $l$ times continuously differentiable
complex-valued functions on $\R^{4n}$ with compact support), and  $g_{\R^{4n}}\in BC^{l}(\R^{4n})$ (the set of bounded continuous 
complex-valued functions on $\R^{4n}$ with all derivatives up to order $l$ bounded and continuous), with a sufficiently large value of $l$.
\end{remark}

Finally, let us discuss what is possible to extract directly from Theorem 2 \cite{coburn:92} and how our theorem relates to that. 

Recall that a fiber of $\pi$ is $\R^{4n}$ with the complex structure $aI+bJ+cK$. In our considerations, $P=(a,b,c)$ is an arbitrary point in ${\mathbb{S}}^2-\{ (a_0,b_0,c_0)\}$, where $(a_0,b_0,c_0)$ is a fixed point in the sphere. There is a Berezin-Toeplitz quantization on each of these $(\R^{4n},aI+bJ+cK)$. Let us explain this in some detail.

As we pointed out before, the map (\ref{diffeo}) is not an $\R$-linear map. However, for each fixed $(a,b,c)$,  
the map (\ref{diffeo}) induces an $\R$-linear isomorphism $S_P:\R^{4n}\to\C^{2n}$, 
${\mathbf{x}}\mapsto ({\mathbf{v}},{\bm{\xi}})$. These are different maps for different $(a,b,c)$.  

The space $\A ^{(k)}(\C^{2n+1})$ is isomorphic to the space $\hat{\A} ^{(k)}(\C^{2n})$ of holomorphic
functions $f$ on $\C^{2n}$ that satisfy
$$
\int_{\C^{2n}}|f({\mathbf{v}},{\bm{\xi}})|^2e^{-k({\mathbf{v}}\cdot\bar{\mathbf{v}}+
  {\bm{\xi}}\cdot\bar{\bm{\xi}})}
   \ 
d\mu({\mathbf{v}},{\bm{\xi}})<\infty .
$$
The isomorphism $\A ^{(k)}(\C^{2n+1})\to\hat{\A }^{(k)}(\C^{2n})$, $h\mapsto \hat{h}$ can be defined by
$$
\hat{h}({\mathbf{v}},{\bm{\xi}})=h({\mathbf{v}},{\bm{\xi}},\zeta)
$$
for a fixed $\zeta\in\C$ 
or by
$$
\hat{h}({\mathbf{v}},{\bm{\xi}})=\int_{\C}h({\mathbf{v}},{\bm{\xi}},\zeta)\frac{1}{\pi(1+|\zeta|^2)^2}d\mu(\zeta).
$$
The inverse map $\hat{h}\mapsto h$ is given by
$$
h({\mathbf{v}},{\bm{\xi}},\zeta)=\hat{h}({\mathbf{v}},{\bm{\xi}})
$$
for all $\zeta \in \C$. 

             Let $k\in\N$. 
             The Toeplitz operator with symbol $h\in L^2(\C^{2n},\Bigl (\frac{k}{\pi}\Bigr )^{2n}e^{-k({\mathbf{v}}\cdot \bar {\mathbf{v}}+{\bm{\xi}}\cdot \bar {\bm{\xi}})}
             d\mu( {\mathbf{v}},{\bm{\xi}}))$ is
             $$
             (T_h^{(k)}f)({\mathbf{v}},{\bm{\xi}})=\Bigl (\frac{k}{\pi}\Bigr )^{2n}\int_{\C^{2n}}h({\mathbf{u}},{\bm{\eta}})f({\mathbf{u}},{\bm{\eta}})e^{k ( {\mathbf{v}}\cdot \bar{\mathbf{u}}+
               {\bm{\xi}}\cdot\bar{\bm{\eta}})}e^{-k({\mathbf{u}}\cdot \bar {\mathbf{u}}+{\bm{\eta}}\cdot \bar {\bm{\eta}})}
             d\mu( {\mathbf{u}},{\bm{\eta}}),
             $$
             Denote, for $l\in \N$, by  $C_c^{l}$ the set of $l$ times continuously differentiable
complex-valued functions on $\C^{2n}$ with compact support, and  by  $BC^{l}$ the set of bounded continuous  
complex-valued functions on $\C^{2n}$ with all derivatives up to order $l$ bounded and continuous. We obtain 
\begin{theorem}\cite[Th. 2]{coburn:92}
  \label{prelimth} Suppose $f$ is in $C_c^{2n+3}$ and $g$ is in $BC^{4n+6}$. Then 
  there is a constant
  $C=C(f,g)$ such that
  $$
  ||T_f^{(k)}T_g^{(k)}-T_{fg}^{(k)}+\frac{1}{k}T^{(k)}_
  {\sum_{j=1}^n\Bigl ( \frac{\partial f}{\partial v_j}\frac{\partial g}{\partial \overline{ v _j}} +
      \frac{\partial f}{\partial \xi_j}\frac{\partial g}{\partial \overline{ \xi_j}}
      \Bigr )}||
  \le \frac{C}{k^2},
  $$
  for all $k>0$.
\end{theorem}
\begin{corollary}
  Suppose $f$ and g are in $C_c^{4n+6}$. Then 
  there is a constant
  $C=C(f,g)$ such that
  $$
  ||ik[T_f^{(k)},T_g^{(k)}]-T^{(k)}_{ \{ f,g\} } ||\le \frac{C}{k}
  $$
    for all $k>0$, where
    $$
  \{ f,g\} = -i\sum_{l=1}^n\Bigr (
\frac{\partial f}{\partial v_l}\frac{\partial g}{\partial \bar v_l} -
 \frac{\partial g}{\partial v_l}\frac{\partial f}{\partial \bar v_l} +
 \frac{\partial f}{\partial \xi_l}\frac{\partial g}{\partial \bar \xi_l}-
 \frac{\partial g}{\partial \xi_l}\frac{\partial f}{\partial \bar \xi_l}
\Bigr ).
$$  
is the standard Poisson bracket on $\C^{2n}$.      
\end{corollary}
The subtle issue is that the functions $f$ and $g$ that appear in Theorem \ref{prelimth} are functions on $\C^{2n}$ and not on $\R^{4n}$. 

Suppose we would like to do quantization on $\R^{4n}$, meaning that to a smooth compactly supported function on $\R^{4n}$ we would like to assign a linear operator, and we would like to use Theorem \ref{prelimth} to ensure that the correspondence principle (Poisson bracket to commutator)  holds.  

If we have two smooth compactly supported functions $f_0$ and $g_0$ on $\R^{4n}$, and we want to use Theorem \ref{prelimth}, then the theorem would be applied to  
the functions $f=f_0\circ (S_P)^{-1}$ and $g=g_0\circ (S_P)^{-1}$. 
\begin{example}
Let $n=1$, $(a_0,b_0,c_0)=(-1,0,0)$. So, we use (\ref{complexcoord1}). Let $f_0:\R^4\to \C$ be defined by $f_0(x_1,x_2,x_3,x_4)=x_1$. 

Take $P=(1,0,0)$. Then $\zeta=0$, $v=z$, $\xi=w$, and the function $f=f_0\circ (S_P)^{-1}$ is defined by $f(v,\xi)=Re(v)$. 

Take, instead, $P=(0,0,-1)$. Then $\zeta=1$, $v=z+\bar w$, $\xi=w-\bar z$, and the function $f=f_0\circ (S_P)^{-1}$ is defined by $f(v,\xi)=\frac{1}{2}
Re(v-\bar\xi)$. 
\end{example}
To recap, suppose we start with a smooth compactly supported function $f_0$ on $\R^{4n}$. Then we get a family of functions $f$, parametrized by  
$P\in {\mathbb{S}}^2-\{ pt\}$ and the corresponding family of Toeplitz operators $T_f^{(k)}$. Thus, given $f_0$ and $k\in\N$, we get "many operators", not one.  
In other words, from \cite{coburn:92} we get a family of Berezin-Toeplitz quantizations, one quantization on each fiber of the twistor projection (with the assumption that we removed one fiber at the beginning of our discussion).  That's a family of Berezin-Toeplitz quantizations, parametrized by $\C$
(one quantization for every complex structure corresponding to a point in $\C$).

The objective of Theorem \ref{mainth} is to get, instead, \underline{one} 
quantization that incorporates all these complex structures at once. Given a smooth compactly supported function $f$ on $\R^{4n}$, we "modify" $f$ (using a diffeomorphism between $\R^{4n+2}$ and $\C^{2n+1}$ and a map from smooth compactly supported functions on $\R^{4n+2}$ to an appropriate function space on $\C^{2n+1}$), to obtain another function, $f_{red}$, and for the operators whose symbols are  these "modified" functions, Theorem \ref{mainth} provide "Berezin-Toeplitz-like" asymptotics. Informally speaking, the total twistor family $\simeq \C^{2n+1}$ and its complex structure includes all copies of the original $\R^{4n}$ with their complex structures. Let $k\in\N$. By assigning, to a function on $\R^{4n}$, an operator on this $\C^{2n+1}$ (the operator $T_{f_{red}}^{(k)}$), we get \underline{one} quantization that captures all the complex structures at once.   

One difference between this approach and the approach in \cite{barron:17}, \cite{castejon:16}, is that in \cite{barron:17}, \cite{castejon:16}, for a given symbol $f$, the quantum operator was built from the three Toeplitz operators (for the three complex structures $I$, $J$, $K$). Here, given $f$, we do have a family of Toeplitz operators  parametrized by $\C$, but we do not use these Toeplitz operators (we do not attempt to build one operator out of all these). Instead, we construct one quantum operator associated to $f$ from different considerations. Also, in \cite{barron:17}, \cite{castejon:16}, it was quantization for \underline{three} complex structures at once. In this paper, we provide a simultaneous quantization for $S^2-\{ pt\}$ of complex structures. As we explained, the point (the single complex structure) removed from $S^2$ can be any point.


\begin{thebibliography}{999}

 \bibitem[AMR]{abraham:83}R. Abraham, J. Marsden, T. Ratiu.
  \newblock {\em Manifolds, tensor analysis, and applications.}
  \newblock Addison-Wesley Publishing Co., Reading, Mass., 1983. 

 \bibitem[AHS]{atiyah:78}M. Atiyah, N. Hitchin, I. Singer. 
  \newblock {\em Self-duality in four-dimensional Riemannian geometry.}
  \newblock Proc. Roy. Soc. London Ser. A 362 (1978), no. 1711, 425-461. 
  

\bibitem[B]{barron:18}T. Barron.
  \newblock {\em Toeplitz operators on K\"ahler manifolds. Examples.}
     \newblock SpringerBriefs in Mathematics. Springer, Cham, 2018. 
  
   \bibitem[BSe]{barron:17}T. Barron, B. Serajelahi.
     \newblock {\em Berezin-Toeplitz quantization, hyperk\"ahler manifolds, and multisymplectic manifolds.} 
     \newblock Glasg. Math. J. 59 (2017), no. 1, 167-187.
     \newblock https://arxiv.org/abs/1406.5159
  
\bibitem[BC]{berger:86}C. Berger, L. Coburn. 
 \newblock {\em Toeplitz operators and quantum mechanics.}
 \newblock J. Funct. Anal. {\bf 68} (1986), no. 3, 273-299. 

\bibitem[CD]{castejon:16} H. Castej\'on D\'iaz.
   \newblock {\em Berezin-Toeplitz quantization on $K3$ surfaces and hyperk\"ahler Berezin-Toeplitz quantization.} 
  \newblock Ph.D. thesis, Universit\'e du Luxembourg, 2016. 
 
\bibitem[Cob]{coburn:92}L. Coburn. 
  \newblock {\em Deformation estimates for the Berezin-Toeplitz quantization.}  
\newblock Commun. Math. Phys. {\bf 149} (1992), 415-424. 

\bibitem[Con]{conway:78}J. Conway.
  \newblock {\em   Functions of one complex variable. Second edition.}
  \newblock Springer-Verlag, New York-Berlin, 1978.

\bibitem[H1]{hitchin:13}N. Hitchin.
\newblock {\em On the hyperk\"ahler/quaternion K\"ahler correspondence.}
\newblock Comm. Math. Phys. {\bf 324} (2013), no. 1, 77-106. 

\bibitem[H2]{hitchin:14}N.Hitchin. 
  \newblock {\em The hyperholomorphic line bundle.}
  \newblock In {\em Algebraic and complex geometry}, pp. 209-223; 
\newblock Springer Proc. Math. Stat., 71, Springer, Cham, 2014. 

\bibitem[Ka]{kaledin:98}D. Kaledin.
  \newblock {\em Integrability of the twistor space for a hypercomplex manifold.} 
  \newblock Selecta Math. (N.S.) {\bf 4} (1998), no. 2, 271-278.
  
\bibitem[KV]{kaled-verbit} D. Kaledin, M. Verbitsky.
  \newblock{\em Non-Hermitian Yang-Mills connections.}
  \newblock Selecta Math. (N.S.) \textbf{4} (1998), 279-320.

\bibitem[Kr]{krantz:92}S. Krantz.
  \newblock {\em Function theory of several complex variables.} Second edition.
  \newblock The Wadsworth and Brooks/Cole Mathematics Series.
   \newblock Pacific Grove, CA, 1992. 
  
\bibitem[MS]{mcduff:98}D. McDuff, D. Salamon. 
  \newblock {\em Introduction to symplectic topology.} Second edition.
  \newblock  Oxford Mathematical Monographs. The Clarendon Press, Oxford University Press, New York, 1998. 

\bibitem[P1]{penrose:68}R. Penrose. 
  \newblock {\em Twistor quantization and curved space-time.}
  \newblock Internat. J. Theor. Phys. Vol. 1, No. 1 (1968), 61-99. 
  
\bibitem[P2]{penrose:99}R. Penrose. 
  \newblock {\em The central programme of twistor theory.}
  \newblock Chaos Solitons Fractals 10 (1999), no. 2-3, 581-611. 
  
\bibitem[R]{range:86}R. Range.
 \newblock {\em Holomorphic functions and integral representations in several complex variables.}
 \newblock  Springer-Verlag, New York, 1986. 
  
\bibitem[S]{salamon:82}S. Salamon.
\newblock {\em Quaternionic K\"ahler manifolds.}
\newblock Invent. Math. {\bf 67} (1982), no. 1, 143-171.

\bibitem[T]{tomberg:15}A. Tomberg.
  \newblock {\em Twistor spaces of hypercomplex manifolds are balanced.} 
  \newblock Adv. Math. {\bf 280} (2015), 282-300.

\bibitem[V1]{verbitsky:95}M. Verbitsky.
\newblock {\em Hyperk\"ahler embeddings and holomorphic symplectic geometry II.}
\newblock GAFA 5 (1995), no. 1, 92-104. 

\bibitem[V2]{verbitsky:15}M. Verbitsky. 
\newblock {\em Degenerate twistor spaces for hyperk\"ahler manifolds.}
\newblock J. Geom. Phys. 91 (2015), 2-11. 

\end{thebibliography}
\end{document}